\documentclass[11pt]{article}%
\usepackage{amsfonts}
\usepackage{amsmath}
\usepackage{amssymb}
\usepackage{amsxtra}
\usepackage{graphicx}
\usepackage{geometry}
\usepackage{color}
\usepackage{colortbl}%
\providecommand{\U}[1]{\protect\rule{.1in}{.1in}}
\newtheorem{theorem}{Theorem}

\newtheorem{corollary}[theorem]{Corollary}

\newtheorem{lemma}[theorem]{Lemma}

\newtheorem{proposition}[theorem]{Proposition}
\newtheorem{remark}[theorem]{Remark}

\numberwithin{equation}{section}
\ifx\pdfoutput\relax\let\pdfoutput=\undefined\fi
\newcount\msipdfoutput
\ifx\pdfoutput\undefined\else
\ifcase\pdfoutput\else
\ifx\paperwidth\undefined\else
\ifdim\paperheight=0pt\relax\else\pdfpageheight\paperheight\fi
\ifdim\paperwidth=0pt\relax\else\pdfpagewidth\paperwidth\fi
\fi\fi\fi
\begin{document}

\title{Thin elastic plates supported over small areas \\II. Variational-asymptotic models.}
\author{G. Buttazzo\thanks{Universit\`{a} di Pisa, Department of Mathematics, Largo B.
Pontecorvo, 5, 56127 Pisa, Italy; email: buttazzo@dm.unipi.it.},
G.Cardone\thanks{Universit\`{a} del Sannio, Department of Engineering, Corso
Garibaldi, 107, 82100 Benevento, Italy; email: giuseppe.cardone@unisannio.it.}%
, S.A.Nazarov\thanks{Mathematics and Mechanics Faculty, St. Petersburg State
University, Universitetskaya nab., 7--9, St. Petersburg, 199034, Russia; Peter
the Great Saint-Petersburg State Polytechnical University, Laboratory
\textquotedblleft Mechanics of New Nano-materials\textquotedblright,
Polytechnicheskaya ul., 29, St. Petersburg, 195251, Russia; Institute of
Problems of Mechanical Engineering RAS, Laboratory \textquotedblleft Mathematical Methods in Mechanics of Materials\textquotedblright, V.O.,
Bolshoj pr., 61, St. Petersburg, 199178, Russia; email: s.nazarov@spbu.ru,
srgnazarov@yahoo.co.uk.}}
\maketitle

\begin{abstract}
\medskip An asymptotic analysis is performed for thin anisotropic elastic
plate clamped along its lateral side and also supported at a small area
$\theta_{h}$ of one base with diameter of the same order as the plate
thickness $h\ll1.$ A three-dimensional boundary layer in the vicinity of the
support $\theta_{h}$ is involved into the asymptotic form which is justified
by means of the previously derived weighted inequality of Korn's type provides
an error estimate with the bound $ch^{1/2}\left\vert \ln h\right\vert .$
Ignoring this boundary layer effect reduces the precision order down to
$\left\vert \ln h\right\vert ^{-1/2}.$ A two-dimensional
variational-asymptotic model of the plate is proposed within the theory of
self-adjoint extensions of differential operators. The only characteristics of
the boundary layer, namely the elastic logarithmic potential matrix of size
$4\times4,$ is involved into the model which however keeps the precision order
$h^{1/2}\left\vert \ln h\right\vert $ in certain norms. Several formulations
and applications of the model are discussed.

\bigskip


\medskip


\end{abstract}

\section{Introduction\label{sect1}}

\subsection{Motivation\label{sect1.1}}

\begin{figure}[ptb]
\begin{center}
\includegraphics[scale=0.65]{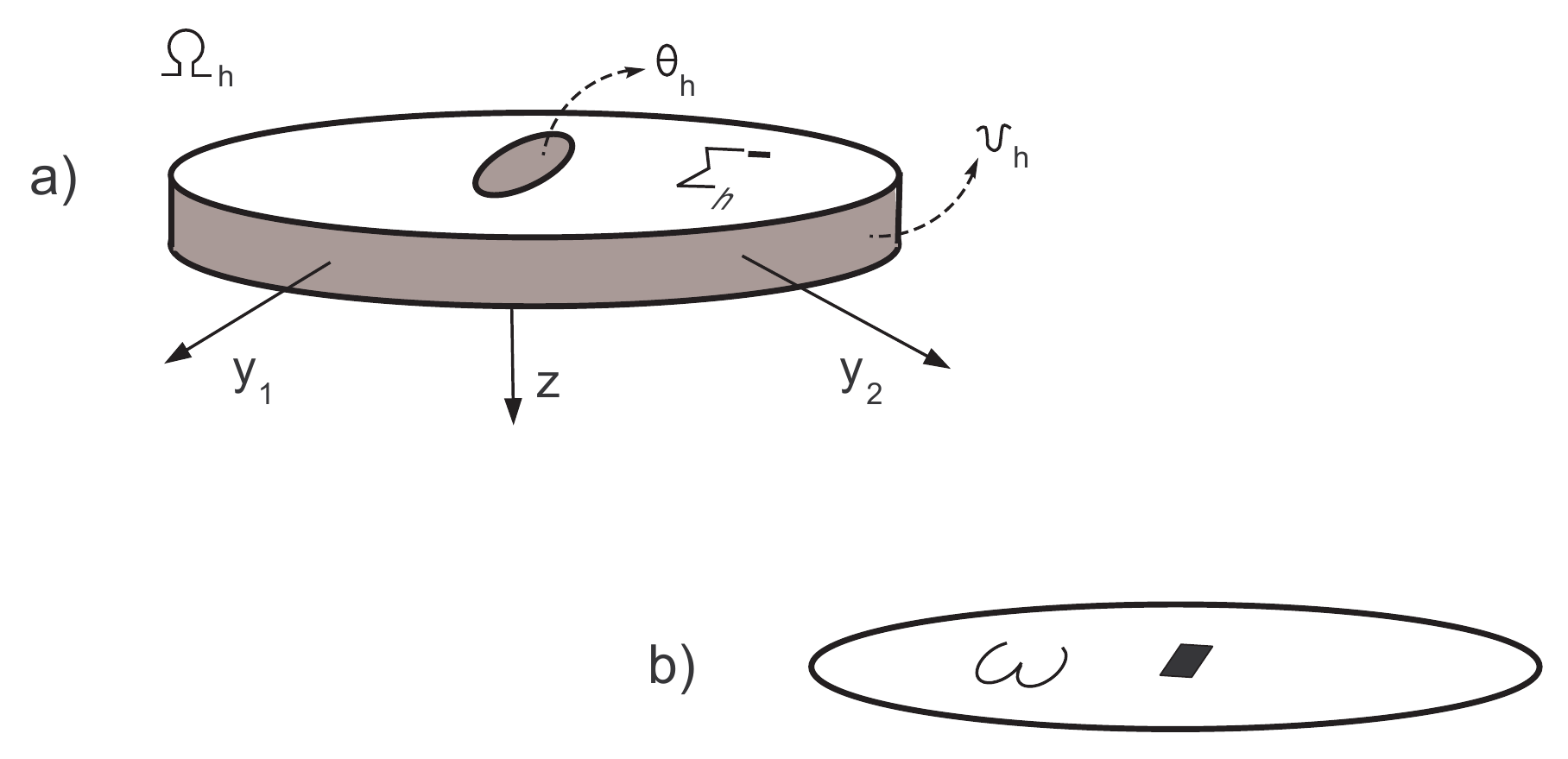}
\end{center}
\caption{A plate with small clamped area (a) and its two-dimensional model
where the black square indicates an asymptotic point condition (b).}%
\label{f1}%
\end{figure}

We consider a thin elastic anisotropic plate
\begin{equation}
\Omega_{h}=\left\{  x=(y,z)\in\mathbb{R}^{2}\times\mathbb{R}:y=\left(
y_{1},y_{2}\right)  \in\omega,\ \left\vert z\right\vert <h/2\right\}
\label{B1.1}%
\end{equation}
of the relative thickness $h\ll1$ with one or several small, of diameter
$O(h),$ support areas at its lower base $\Sigma_{h}^{-}$ (see Figure
\ref{f1},a). In this paper we assume that the plate is clamped along the
lateral side $\upsilon_{h}$ while the bases $\Sigma_{h}^{\pm}$ are
traction-free, except at the only small support zone%
\begin{equation}
\theta_{h}=\left\{  x:\eta:=h^{-1}y\in\theta,\ z=-h/2\right\}  . \label{B1.3}%
\end{equation}
Here, $\omega$ and $\theta$ are domains in the plane $\mathbb{R}^{2}$ with
smooth boundaries and compact closures (see Figure \ref{f1},b), the
$y$-coordinate origin $\mathcal{O}$ is put inside $\theta$,
\begin{equation}
\Sigma_{h}^{\pm}=\left\{  x:y\in\omega,\ z=\pm h/2\right\}  ,\ \ \ \upsilon
_{h}=\left\{  x:y\in\partial\omega,\ \left\vert z\right\vert <h/2\right\}  .
\label{B1.2}%
\end{equation}

This paper is a direct sequel of \cite{ButCaNa1} where a convenient form of
Korn's inequality in $\Omega_{h}$ was derived as well as a boundary layer
effect near $\theta_{h}$ was investigated scrupulously. These results,
combined with the classical procedure of dimension reduction, see
\cite{Ciarlet1, Ciarlet2, LeDret, Nabook, SaPa} etc., and an asymptotic
analysis of singular perturbation of boundaries, see \cite{Ilin, MaNaPl},
allow us here to construct and justify asymptotic expansions of elastic fields
in $\Omega_{h},$ i.e., displacements, strains and stresses.

A small Dirichlet perturbation of boundary conditions may lead to neglectable
changes in a solution but in our case, unfortunately, such a change is fair
and significant in two aspects. First of all, fixing the plate (\ref{B1.1})
with the sharp supports at $\theta_{h}^{1},...,\theta_{h}^{J},$ cf.
(\ref{B1.3}), results in the Sobolev point conditions, see, e.g., \cite{BuNa},
namely the average deflexion $w_{3}$ vanishes at the points $\mathcal{O}%
^{1},...,\mathcal{O}^{J},$ the centers of the small support areas. Due to the
Sobolev embedding theorem in $\mathbb{R}^{2}$ these additional restrictions
are natural in the variational formulation of the traditional two-dimensional
model of Kirchhoff's plate while the only implication defect turns into a
certain lack of smoothness properties of the solution at the points
$\mathcal{O}^{1},...,\mathcal{O}^{J}.$ Much more disagreeable input of the
small support areas becomes a crucial reduction of the convergence rate:
instead of the the conventional rate $O(\sqrt{h})$ for the Kirchhoff model, we
detect $O(|\ln h|^{-1/2}),$ cf. \cite{na472}, that is unacceptable for any
application purpose.

The three-dimensional boundary layers, described and investigated in
\cite{ButCaNa1}, see also \cite{na135, na243}, help to construct specific
corerctions in the vicinity of $\theta_{h},$ to take into account their
influence on the far-field asymptotics and, as a result, to achieve the
appropriate error estimate $O(\sqrt{h}|\ln h|).$ At the same time, the global
asymptotic approximation of the elastic fields becomes quite complicated so
that its utility in application is rather doubtful. In this way, it is worth
to restrict our analysis to the fair-field asymptotics whose dependence on the
small parameter $h$ remains complicated but can be expressed through rational
functions in $\ln h$ and an integral characteristics of the support area,
namely a $4\times4$-matrix of elastic capacity \cite{ButCaNa1, na135}.

To provide a correct and mathematically rigorous two-dimensional model of the
supported plate, we employ two approaches. First, we turn to the technique of
self-adjoint extensions of differential operators (see the pioneering paper
\cite{BeFa}, the review \cite{Pav} and, e.g., the papers \cite{KaPa, na101,
na341, na366}, especially \cite{na239} connecting this approach with the
method of matched asymptotic expansions). To simulate the influence of
localized special boundary layer, we provide a proper choice of singular
solutions by selecting a special self-adjoint extension in the Lebesgue space
$L^{2}(\omega)$ of the matrix $\mathcal{L}(\nabla_{y})$ of differential
operators in the Kirchhoff plate model, see (\ref{3.28}), (\ref{3.27}).
Parameters of this extension depend on the quantity $|\ln h|$ and the elastic
logarithmic capacity matrix, see \cite[Sect. 3.2]{ButCaNa1}, which is similar
to the logarithmic capacity in harmonic analysis, cf. \cite{Land, PoSe}.

The introduced model is two-fold and can be formulated as either an abstract
equation with the obtained self-adjoint operator, or a problem on a stationary
point of the natural energy functional in a weighted function space with
detached asymptotics. Advantages and deficiencies of these approaches are
discussed in Section \ref{sect4}. In any case the compelled removal of a major
part of the three-dimensional layer from the asymptotic form surely reduces
capabilities of the model which, in particular, does not give the global
approximation of stresses and strains, but does outside a neighborhood of the
sets $\overline{\theta_{h}}$ and $\overline{\upsilon_{h}}$ where the Dirichlet
conditions are imposed. Thereupon, we mention that, first, the Kirchhoff model
itself does not approximate correctly stresses and strains near the clamped
edge, too, and, second, quite many applications demand to determine
asymptotics of displacements only. In Corollaries \ref{cor5.99} and
\ref{cor5.77} we will give some error estimates for the model, in particular,
an asymptotic formula for the elastic energy of the three-dimensional
thin-plate. Results \cite{na400} predict that the model also can furnish
asymptotics of eigenfrequencies of the plate but this will be a topic of an
oncoming paper.

\subsection{Formulation of the problem\label{sect1.2}}

Similarly to \cite{ButCaNa1}, in the domain (\ref{B1.1}) we consider the mixed
boundary-value problem of the elasticity theory%
\begin{align}
L(\nabla)u(h,x)  &  :=D\left(  -\nabla\right)  ^{\top}AD\left(  \nabla\right)
u\left(  h,x\right)  =f\left(  h,x\right)  ,\ \ \ x\in\Omega_{h}%
,\label{B1.11}\\
N^{+}(\nabla)u(h,x)  &  :=D\left(  e_{3}\right)  ^{\top}AD\left(
\nabla\right)  u\left(  h,x\right)  =0,\ \ \ x\in\Sigma_{h}^{+},
\label{B1.12}\\
N^{-}(\nabla)u(h,x)  &  :=D\left(  -e_{3}\right)  ^{\top}AD\left(
\nabla\right)  u\left(  h,x\right)  =0,\ \ \ x\in\Sigma_{h}^{\bullet}%
=\Sigma_{h}^{-}\setminus\theta_{h},\nonumber\\
u\left(  h,x\right)   &  =0,\ \ \ x\in\theta_{h},\label{B1.13}\\
u\left(  h,x\right)   &  =0,\ \ \ x\in\upsilon_{h}. \label{B1.14}%
\end{align}
Let us explain the Voigt-Mandel notation we use throughout the paper. The
displacement vector $\left(  u_{1}(x),u_{2}(x),u_{3}(x)\right)  ^{\top}$ is
regarded as a column in $\mathbb{R}^{3}$, $\top$ stands for transposition and
the Hooke's law
\begin{equation}
\sigma\left(  u\right)  =A\varepsilon\left(  u\right)  =AD\left(
\nabla\right)  u \label{B1.9}%
\end{equation}
relates the strain column%
\[
\varepsilon=\left(  \varepsilon_{11},\varepsilon_{22},2^{1/2}\varepsilon
_{12},2^{1/2}\varepsilon_{13},2^{1/2}\varepsilon_{23},\varepsilon_{33}\right)
^{\top},\ \ \ \varepsilon_{jk}\left(  u\right)  =\frac{1}{2}\left(
\frac{\partial u_{j}}{\partial x_{k}}+\frac{\partial u_{k}}{\partial x_{j}%
}\right)  ,
\]
with the stress column $\sigma\left(  u\right)  $ of the same structure as
$\varepsilon\left(  u\right)  .$ Moreover, the stiffness matrix $A$ of size
$6\times6,$ symmetric and positive definite, figures in (\ref{B1.9}) as well
as the differential operators%
\begin{equation}
D\left(  \nabla\right)  ^{\top}=\left(
\begin{array}
[c]{cccccc}%
\partial_{1} & 0 & 2^{-1/2}\partial_{2} & 2^{-1/2}\partial_{3} & 0 & 0\\
0 & \partial_{2} & 2^{-1/2}\partial_{1} & 0 & 2^{-1/2}\partial_{3} & 0\\
0 & 0 & 0 & 2^{-1/2}\partial_{1} & 2^{-1/2}\partial_{2} & \partial_{3}%
\end{array}
\right)  ,\ \ \ \partial_{j}=\frac{\partial}{\partial x_{j}},\ \ \ \nabla
=\left(
\begin{array}
[c]{c}%
\partial_{1}\\
\partial_{2}\\
\partial_{3}%
\end{array}
\right)  . \label{B1.7}%
\end{equation}
The variational formulation of problem (\ref{B1.11})-(\ref{B1.14}) implies the
integral identity
\begin{equation}
\left(  AD\left(  \nabla\right)  u,D\left(  \nabla\right)  v\right)
_{\Omega_{h}}=\left(  f,v\right)  _{\Omega_{h}},\ \ \ \forall v\in H_{0}%
^{1}\left(  \Omega_{h};\Gamma_{h}\right)  ^{3}, \label{B1.16}%
\end{equation}
where $\left(  \ ,\ \right)  _{\Omega_{h}}$ is the scalar product in
$L^{2}\left(  \Omega_{h}\right)  ,$ $H_{0}^{1}\left(  \Omega_{h};\Gamma
_{h}\right)  $ is a subspace of functions in the Sobolev space $H^{1}\left(
\Omega_{h}\right)  $ which vanish on $\Gamma_{h}=\theta_{h}\cup\upsilon_{h},$
see (\ref{B1.13}) and (\ref{B1.14}), and the last superscript $3$ in
(\ref{B1.16}) indicates the number of components in the test (vector) function
$v=\left(  v_{1},v_{2},v_{3}\right)  ^{\top}.$

We here have listed only notation of further use while more explanation and an
example of the isotropic material can be found in Section 1.2 of
\cite{ButCaNa1}.

\subsection{Architecture of the paper\label{sect1.5}}

In \S \ref{sect3} we present the standard asymptotic procedure of dimension
reduction to provide explicit formulas in the Kirchhoff plate's model. Based
on the results obtained in \cite{ButCaNa1}, we construct in \S \ref{sect4}
asymptotic expansions of a solution $u(h,x)$ to problem (\ref{B1.11}%
)-(\ref{B1.14}) including the boundary layer near the support area
(\ref{B1.3}) and derive the error estimates in various versions.

The technique of self-adjoint extensions of differential operators is outlined
in \S \ref{sect4} for the Kirchhoff plate supported at one point, i.e. with
the Sobolev condition, together with the choice of the extension parameters on
the base of the previous asymptotic analysis. Finally, the
asymptotic-variational model mentioned in Section \ref{sect1.1} is discussed
at length in \S \ref{sect5}.

\section{The limit two-dimensional problem\label{sect3}}

\subsection{The dimension reduction\label{sect3.1}}

We suppose that the right-hand side of system (\ref{B1.11}) takes the form%
\begin{equation}
f\left(  h,x\right)  =h^{-1/2}f^{0}\left(  y,\zeta\right)  +h^{1/2}%
f^{1}\left(  y\right)  +\widetilde{f}\left(  h,x\right)  \label{3.1}%
\end{equation}
where $\zeta=h^{-1}z\in\left(  -1/2,1/2\right)  $ is the stretched transverse
coordinate and $f^{0},$ $f^{1}$ are vector functions subject to the conditions%
\begin{equation}
\int_{-1/2}^{1/2}f_{3}^{0}\left(  y,\zeta\right)  d\zeta=0,\ \ \ f_{1}%
^{1}\left(  y\right)  =f_{2}^{1}\left(  y\right)  =0,\ \ \ y\in\omega.
\label{3.2}%
\end{equation}
Now we assume $f^{0}$ and $f^{1}$ smooth and postpone a detailed description
of properties of $f^{0}$, $f^{1}$ and $\widetilde{f}$ to Section
\ref{sect4.2}. We emphasize that, owing to linearity of the problem, any
reasonable right-hand side can be represented as in (\ref{3.1}), (\ref{3.2}).
After appropriate rescaling our choice of the factors $h^{-1/2}$ and $h^{1/2}$
in (\ref{3.1}) as well as the orthogonality condition in (\ref{3.2}) ensures
that the potential energy stored by the plate $\Omega_{h}$ under the volume
force (\ref{3.1}) gets order $1=h^{0},$ see Section \ref{sect5.4}.

As shown in \cite{na273}, \cite[Ch.4]{Nabook}, the asymptotic ansatz
\begin{equation}
u\left(  h,x\right)  =h^{-3/2}U^{0}\left(  y\right)  +h^{-1/2}U^{1}\left(
y,\zeta\right)  +h^{1/2}U^{2}\left(  y,\zeta\right)  +... \label{3.3}%
\end{equation}
of the solution to problem (\ref{B1.11})-(\ref{B1.13}) is perfectly adjusted
with decomposition (\ref{3.1}). The coordinate change $x\mapsto\left(
y,\zeta\right)  $ splits the differential operators $L\left(  \nabla\right)  $
and $N^{\pm}\left(  \nabla\right)  $ in (\ref{B1.11}) and (\ref{B1.12}),
respectively, as follows:%
\begin{align}
L\left(  \nabla\right)   &  :=D\left(  -\nabla\right)  ^{\top}AD\left(
\nabla\right)  =h^{-2}L^{0}\left(  \partial_{\zeta}\right)  +h^{-1}%
L^{1}\left(  \nabla_{y},\partial_{\zeta}\right)  +h^{0}L^{2}\left(  \nabla
_{y}\right)  ,\label{3.4}\\
N^{\pm}\left(  \nabla\right)   &  :=D\left(  0,0,\pm1\right)  ^{\top}AD\left(
\nabla\right)  =h^{-1}N^{0\pm}\left(  \partial_{\zeta}\right)  +h^{0}N^{1\pm
}\left(  \nabla_{y}\right)  ,\nonumber
\end{align}
where%
\begin{align}
L^{0}  &  =-D_{\zeta}^{\top}AD_{\zeta},\ \ \ L^{1}=-D_{\zeta}^{\top}%
AD_{y}-D_{y}^{\top}AD_{\zeta},\ \ \ L^{2}=-D_{y}^{\top}AD_{y},\nonumber\\
N^{0\pm}  &  =\pm D_{3}^{\top}AD_{\zeta},\ \ \ N^{1\pm}=\pm D_{3}^{\top}%
AD_{y}\nonumber\\
D_{\zeta}  &  =D\left(  0,0,\partial_{\zeta}\right)  ,\ \ \ D_{y}=D\left(
\nabla_{y},0\right)  ,\ \ \ D_{3}=\left(  0,0,1\right)  ,\label{3.6}\\
\partial_{\zeta}  &  =\partial/\partial\zeta,\ \ \ \nabla_{y}=\left(
\partial/\partial y_{1},\partial/\partial y_{2}\right)  .\nonumber
\end{align}
We now insert formulas (\ref{3.4}) and (\ref{3.3}), (\ref{3.1}) into
(\ref{B1.11}) and (\ref{B1.12}) and collect coefficients of $h^{p-2}$ and
$h^{p-1},$ respectively. Equalizing their sums yields a recursive sequence of
the Neumann problems on the interval $\left(  -1/2,1/2\right)  \ni\zeta$ with
the parameter $y\in\omega$. For $p=0,1,2,$ we have
\begin{align}
L^{0}U^{p}  &  =F^{p}:=-L^{1}U^{p-1}-L^{2}U^{p-2}\text{ for }\zeta\in\left(
-\tfrac{1}{2},\tfrac{1}{2}\right)  ,\label{3.7}\\
N^{0\pm}U^{p}  &  =G^{p\pm}:=-N^{1\pm}U^{p-1}\text{ at }\zeta=\pm\tfrac{1}%
{2},\nonumber
\end{align}
where $U^{-2}=U^{-1}=0.$ Since $F^{0}=0,\ G^{0\pm}=0$ and $F^{1}=0,\ G^{1\pm
}=\mp D_{3}^{\top}AD_{y}U^{0},$ we readily take%
\begin{equation}
U^{0}\left(  y\right)  =w_{3}\left(  y\right)  e_{3},\ \ \ U^{1}\left(
y\right)  =\sum\limits_{i=1}^{2}e_{i}\left(  w_{i}\left(  y\right)
-\zeta\frac{\partial w_{3}}{\partial y_{i}}\left(  y\right)  \right)
\label{3.8}%
\end{equation}
as solutions to problem (\ref{3.7}) with $p=0$ and $p=1.$ Here, $e_{j}$ is the
unit vector of the $x_{j}$-axis and%
\begin{equation}
w=\left(  w_{1},w_{2},w_{3}\right)  ^{\top} \label{3.9}%
\end{equation}
is a vector function to be determined. Notice that $h^{-3/2}w_{3}\left(
y\right)  $ and $h^{-1/2}w_{i}\left(  y\right)  $ will be the averaged
deflection and longitudinal displacements in the plate.

If $p=2,$ one may easily verify that the right-hand sides in (\ref{3.7})
become%
\begin{align}
F^{2}\left(  y,\zeta\right)   &  =D_{\zeta}^{\top}AD_{y}U^{1}+D_{y}^{\top
}A\left(  D_{\zeta}U^{1}+D_{y}U^{0}\right)  =D_{\zeta}^{\top}A\mathcal{Y}%
\left(  \zeta\right)  \mathcal{D}\left(  \nabla_{y}\right)  w\left(  y\right)
,\label{3.77}\\
G^{2\pm}\left(  y\right)   &  =\mp D_{3}^{\top}A\mathcal{Y}\left(  1/2\right)
\mathcal{D}\left(  \nabla_{y}\right)  w\left(  y\right) \nonumber
\end{align}
where $D_{\zeta}U^{1}+D_{y}U^{0}=0$ in view of (\ref{3.8}) and%
\begin{equation}
\mathcal{Y}\left(  \zeta\right)  =\left(
\begin{array}
[c]{cc}%
\mathbb{I}_{3} & -2^{1/2}\mathbb{I}_{3}\zeta\\
\mathbb{O}_{3} & \mathbb{O}_{3}%
\end{array}
\right)  ,\ \mathcal{D}\left(  \nabla_{y}\right)  ^{\top}=\left(
\begin{array}
[c]{cccccc}%
\partial_{1} & 0 & 2^{-1/2}\partial_{2} & 0 & 0 & 0\\
0 & \partial_{2} & 2^{-1/2}\partial_{1} & 0 & 0 & 0\\
0 & 0 & 0 & 2^{-1/2}\partial_{1}^{2} & 2^{-1/2}\partial_{2}^{2} & \partial
_{1}\partial_{2}%
\end{array}
\right)  , \label{3.10}%
\end{equation}
cf. (\ref{B1.7}), $\mathbb{I}_{3}$ and $\mathbb{O}_{3}$ are the unit and null
matrices of size $3\times3.$

Hence,%
\begin{equation}
U^{2}\left(  y,\zeta\right)  =\mathcal{X}\left(  \zeta\right)  \mathcal{D}%
\left(  \nabla_{y}\right)  w\left(  y\right)  \label{3.11}%
\end{equation}
and $\mathcal{X}$ is a matrix solution (of size $3\times6)$ to the problem%
\begin{equation}
-D_{\zeta}^{\top}AD_{\zeta}\mathcal{X}\left(  \zeta\right)  =D_{\zeta}^{\top
}A\mathcal{Y}\left(  \zeta\right)  ,\ \zeta\in\left(  -\tfrac{1}{2},\tfrac
{1}{2}\right)  ,\ \pm D_{3}^{\top}AD_{\zeta}\mathcal{X}\left(  \pm\tfrac{1}%
{2}\right)  =\mp D_{3}^{\top}A\mathcal{Y}\left(  \pm\tfrac{1}{2}\right)  .
\label{3.12}%
\end{equation}
Compatibility conditions in problem (\ref{3.12}) are evidently satisfied.

The problem with $p=3$%
\begin{equation}
L^{0}U^{p}=F^{p}\text{ for }\zeta\in\left(  -\tfrac{1}{2},\tfrac{1}{2}\right)
,\text{\ \ \ }N^{0\pm}U^{p}=G^{p\pm}\text{ at }\zeta=\pm\tfrac{1}{2}
\label{3.99}%
\end{equation}
gets the right-hand sides%
\begin{align}
F^{3}  &  =-L^{1}U^{2}-L^{0}U^{1}+f^{0}=D_{\zeta}^{\top}AD_{y}U^{2}%
+D_{y}^{\top}A\left(  D_{\zeta}U^{2}+D_{y}U^{1}\right)  +f^{0},\label{3.98}\\
G^{3\pm}  &  =\mp N^{1\pm}U^{2}=\mp D_{3}^{\top}AD_{y}U^{2}.\nonumber
\end{align}

Let us consider the compatibility conditions%
\begin{equation}
I_{j}^{p}\left(  y\right)  :=\int_{-1/2}^{1/2}F_{j}^{p}\left(  y,\zeta\right)
d\zeta+G_{j}^{p+}\left(  y\right)  +G_{j}^{p-}\left(  y\right)
=0,\ \ \ j=1,2,3, \label{3.13}%
\end{equation}
with $p=3.$ By (\ref{3.98}) and (\ref{3.8}), (\ref{3.11}), we have%
\begin{equation}
D_{\zeta}U^{2}\left(  y,\zeta\right)  +D_{y}U^{1}\left(  y,\zeta\right)
=\left(  D_{\zeta}\mathcal{X}\left(  \zeta\right)  +\mathcal{Y}\left(
\zeta\right)  \right)  \mathcal{D}\left(  \nabla_{y}\right)  w\left(
y\right)  \label{3.14}%
\end{equation}
and%
\[
I_{i}^{3}\left(  y\right)  =e_{i}^{\top}D_{y}^{\top}\int_{-1/2}^{1/2}A\left(
D_{\zeta}\mathcal{X}\left(  \zeta\right)  +\mathcal{Y}\left(  \zeta\right)
\right)  d\zeta\mathcal{D}\left(  \nabla_{y}\right)  w\left(  y\right)
d\zeta+\int_{-1/2}^{1/2}f_{i}^{0}\left(  y,\zeta\right)  d\zeta.
\]
Since $e_{i}^{\top}D_{y}^{\top}=e_{i}^{\top}\mathcal{D}\left(  \nabla
_{y}\right)  $ according to (\ref{B1.7}), (\ref{3.6}), (\ref{3.10}), the
equalities $I_{i}^{3}\left(  y\right)  =0$ with $i=1,2$ are nothing but two
upper lines in the system of differential equations%
\begin{equation}
\mathcal{D}\left(  -\nabla_{y}\right)  ^{\top}\mathcal{AD}\left(  \nabla
_{y}\right)  w\left(  y\right)  =g\left(  y\right)  ,\ \ \ y\in\omega,
\label{3.15}%
\end{equation}
where%
\begin{equation}
\mathcal{A}=\int_{-1/2}^{1/2}\mathcal{Y}\left(  \zeta\right)  ^{\top}A\left(
D_{\zeta}\mathcal{X}\left(  \zeta\right)  +\mathcal{Y}\left(  \zeta\right)
\right)  d\zeta\label{3.16}%
\end{equation}
and%
\begin{equation}
g_{i}\left(  y\right)  =\int_{-1/2}^{1/2}f_{i}^{0}\left(  y,\zeta\right)
d\zeta,\ \ \ i=1,2. \label{3.17}%
\end{equation}
The compatibility condition (\ref{3.13}) with $j=3$ is fulfilled
automatically. Indeed, we take into account the orthogonality condition in
(\ref{3.2}) and write%
\begin{align*}
I_{3}^{3}\left(  y\right)   &  =\int_{-1/2}^{1/2}e_{3}^{\top}D_{y}^{\top
}A\left(  D_{\zeta}\mathcal{X}\left(  \zeta\right)  +\mathcal{Y}\left(
\zeta\right)  \right)  d\zeta~\mathcal{D}\left(  \nabla_{y}\right)  w\left(
y\right) \\
&  =\sum\limits_{i=1}^{2}\frac{\partial}{\partial y_{i}}\int_{-1/2}%
^{1/2}\left(  D_{\zeta}\zeta e_{i}\right)  ^{\top}A\left(  D_{\zeta
}\mathcal{X}\left(  \zeta\right)  +\mathcal{Y}\left(  \zeta\right)  \right)
d\zeta~\mathcal{D}\left(  \nabla_{y}\right)  w\left(  y\right) \\
&  =\sum\limits_{i=1}^{2}\frac{\partial}{\partial y_{i}}\Bigg(-\int
_{-1/2}^{1/2}e_{i}^{\top}\zeta D_{\zeta}^{\top}A\left(  D_{\zeta}%
\mathcal{X}\left(  \zeta\right)  +\mathcal{Y}\left(  \zeta\right)  \right)
d\zeta\\
&  \quad+e_{i}^{\top}\zeta D_{3}^{\top}A\left(  D_{\zeta}\mathcal{X}\left(
\zeta\right)  +\mathcal{Y}\left(  \zeta\right)  \right)  \left\vert
_{\zeta=-1/2}^{1/2}\right.  \Bigg)\mathcal{D}\left(  \nabla_{y}\right)
w\left(  y\right)  =0.
\end{align*}
Here, we again applied representation (\ref{3.14}) together with the identity%
\begin{equation}
D_{y}e_{3}=\sum\limits_{i=1}^{2}D_{\zeta}\zeta e_{i}\frac{\partial}{\partial
y_{i}} \label{3.19}%
\end{equation}
inherited from (\ref{B1.7}), (\ref{3.6}), and finally recalled problem
(\ref{3.12}) for the $3\times6$-matrix $\mathcal{X}$.

Let us demonstrate that the third line of system (\ref{3.15}) is equivalent to
the equality%
\begin{equation}
I_{3}^{4}\left(  y\right)  =0, \label{3.20}%
\end{equation}
that is the third compatibility condition (\ref{3.13}) in problem (\ref{3.99})
with the right-hand sides%
\begin{align*}
F^{3}  &  =-L^{1}U^{3}-L^{0}U^{2}+f^{1}=D_{\zeta}^{\top}AD_{y}U^{3}%
+D_{y}^{\top}A\left(  D_{\zeta}U^{3}+D_{y}U^{2}\right)  +f_{3}^{1},\\
G^{4\pm}  &  =\mp N^{1\pm}U^{3}=\mp D_{3}^{\top}AD_{y}U^{3}.
\end{align*}
We have%
\begin{align}
I_{3}^{4}\left(  y\right)   &  =\int_{-1/2}^{1/2}e_{3}^{\top}D_{y}^{\top
}A\left(  D_{\zeta}U^{3}+D_{y}U^{2}\right)  d\zeta+f_{3}^{1}\label{3.21}\\
&  =\sum\limits_{i=1}^{2}\frac{\partial}{\partial y_{i}}\int_{-1/2}%
^{1/2}\left(  D_{\zeta}\zeta e_{i}\right)  ^{\top}A\left(  D_{\zeta}%
U^{3}+D_{y}U^{2}\right)  d\zeta+f_{3}^{1}\nonumber\\
&  =f_{3}^{1}+\sum\limits_{i=1}^{2}\frac{\partial}{\partial y_{i}}e_{i}^{\top
}\left(  -\int_{-1/2}^{1/2}\zeta D_{\zeta}^{\top}A\left(  D_{\zeta}U^{3}%
+D_{y}U^{2}\right)  d\zeta+\left(  \zeta D_{3}^{\top}A\left(  D_{\zeta}%
U^{3}+D_{y}U^{2}\right)  \right)  |_{\zeta=-1/2}^{1/2}\right) \nonumber\\
&  =f_{3}^{1}+\sum\limits_{i=1}^{2}\frac{\partial}{\partial y_{i}}e_{i}^{\top
}\int_{-1/2}^{1/2}\zeta\left(  D_{y}^{\top}A\left(  D_{\zeta}U^{2}+D_{y}%
U^{1}\right)  +f^{0}\right)  d\zeta\nonumber\\
&  =g_{3}+\sum\limits_{i=1}^{2}\frac{\partial}{\partial y_{i}}e_{i}^{\top
}D_{y}^{\top}\int_{-1/2}^{1/2}\zeta A\left(  D_{\zeta}\mathcal{X}%
+\mathcal{Y}\right)  d\zeta\mathcal{D}\left(  \nabla_{y}\right)  w\nonumber
\end{align}
where%
\begin{equation}
g_{3}\left(  y\right)  =f_{3}^{1}\left(  y\right)  +\sum\limits_{i=1}^{2}%
\int_{-1/2}^{1/2}\zeta\frac{\partial f_{i}^{0}}{\partial y_{i}}\left(
y,\zeta\right)  d\zeta. \label{3.22}%
\end{equation}
Let us clarify calculation (\ref{3.21}). First, we used identity (\ref{3.19})
and integrate by parts in the interval $\left(  -1/2,1/2\right)  .$ Then we
took into account\ problem (\ref{3.7}) with $p=3$ and its right-hand sides
(\ref{3.98}). Finally, formula (\ref{3.14}) was applied.

In view of (\ref{B1.7}) and (\ref{3.10}) one observes that%
\[
\sum\limits_{i=1}^{2}\frac{\partial}{\partial y_{i}}e_{i}^{\top}D_{y}^{\top
}\zeta=e_{3}^{\top}\mathcal{D}\left(  \nabla_{y}\right)  \mathcal{Y}\left(
\zeta\right)
\]
and, hence, equality (\ref{3.20}) indeed implies the third line of system
(\ref{3.15}) due to definition (\ref{3.16}).

\subsection{The plate equations\label{sect3.2}}

Owing to the Dirichlet condition (\ref{B1.14}) on the lateral side
$\upsilon_{h}$ of the plate $\Omega_{h},$ we supply system (\ref{3.15}) with
the boundary conditions%
\begin{equation}
w_{i}\left(  y\right)  =0,\ i=1,2,\ \ \ w_{3}\left(  y\right)
=0,\ \ \ \partial_{n}w_{3}\left(  y\right)  =0,\ y\in\partial\omega,
\label{3.23}%
\end{equation}
where $\partial_{n}=n^{\top}\nabla_{y}$ and $n=(n_{1},n_{2})^{\top}$ is the
unit vector of the outward normal at the boundary of the domain $\omega
\subset\mathbb{R}^{2}.$ Notice that conditions (\ref{3.23}) make the main
asymptotic terms (\ref{3.8}) in ansatz (\ref{3.3}) vanish at the lateral
boundary $\upsilon_{h}.$

It follows from the weighted Korn's inequality (2.10) in \cite{ButCaNa1} and
will be shown by other means later that the Dirichlet conditions (\ref{B1.13})
at the small support zones (\ref{B1.3}) lead to the Sobolev (point) condition
at the $y$-coordinates origin $\mathcal{O},$ namely%
\begin{equation}
w_{3}\left(  \mathcal{O}\right)  =0. \label{3.24}%
\end{equation}
To prove the unique solvability of the limit two-dimensional problem
(\ref{3.15}), (\ref{3.23}), (\ref{3.24}) needs a piece of information on the
matrix $\mathcal{A}$ of the system.

\begin{lemma}
\label{lem3.1}Matrix (\ref{3.16}) is block-diagonal and takes the form%
\begin{equation}
\mathcal{A}=\left(
\begin{array}
[c]{cc}%
\mathcal{A}^{\prime} & \mathbb{O}_{3}\\
\mathbb{O}_{3} & \mathcal{A}_{\left(  3\right)  }%
\end{array}
\right)  =\left(
\begin{array}
[c]{cc}%
A^{0} & \mathbb{O}_{3}\\
\mathbb{O}_{3} & \frac{1}{6}A^{0}%
\end{array}
\right)  \label{3.25}%
\end{equation}
where $A^{0}=A_{\left(  yy\right)  }-A_{\left(  yz\right)  }A_{\left(
zz\right)  }^{-1}A_{\left(  zy\right)  }$ is a symmetric and positive definite
$3\times3$-matrix constructed from submatrices of size $3\times3$ in the
representation of the stiffness matrix%
\begin{equation}
A=\left(
\begin{array}
[c]{cc}%
A_{\left(  yy\right)  } & A_{\left(  yz\right)  }\\
A_{\left(  zy\right)  } & A_{\left(  zz\right)  }%
\end{array}
\right)  . \label{3.26}%
\end{equation}

\end{lemma}

\textbf{Proof.} We introduce the diagonal matrix $\mathbb{J}%
=\operatorname*{diag}\left(  2^{-1/2},2^{-1/2},1\right)  $ and, by
(\ref{3.26}) and (\ref{B1.7}), rewrite problem (\ref{3.12}) as follows:%
\begin{align*}
-A_{\left(  zz\right)  }\mathbb{J}\partial_{\zeta}^{2}\mathcal{X}\left(
\zeta\right)   &  =A_{\left(  zy\right)  }(\mathbb{O}_{3},-2^{1/2}%
\mathbb{I}_{3}),\ \zeta\in\left(  -1/2,1/2\right)  ,\ \ \\
\pm A_{\left(  zz\right)  }\mathbb{J}\partial_{\zeta}\mathcal{X}\left(
\pm1/2\right)   &  =\mp A_{\left(  zy\right)  }(\mathbb{I}_{3},\mp
2^{1/2}\mathbb{I}_{3}).
\end{align*}
Its solution thus reads%
\begin{equation}
\mathcal{X}\left(  \zeta\right)  =\mathbb{J}^{-1}A_{\left(  zz\right)  }%
^{-1}A_{\left(  zy\right)  }\left(  -\zeta\mathbb{I}_{3},2^{1/2}\left(
\tfrac{\zeta^{2}}{2}-\tfrac{1}{24}\right)  \mathbb{I}_{3}\right)  .
\label{3.2678}%
\end{equation}
Inserting this formula into (\ref{3.16}) assures representation (\ref{3.25}).
Since matrix (\ref{3.26}) is symmetric and positive definite, these properties
are evidently passed to $A^{0}$. $\blacksquare$

Representation (\ref{3.25}) indicates a remarkable fact: for any homogeneous
anisotropic plate the limit two-dimensional system (\ref{3.15}) divides into
the fourth-order differential operator%
\begin{equation}
\mathcal{L}_{3}\left(  \nabla_{y}\right)  =\tfrac{1}{6}\mathcal{D}_{3}\left(
\nabla_{y}\right)  ^{\top}A^{0}\mathcal{D}_{3}\left(  \nabla_{y}\right)
\label{3.27}%
\end{equation}
for the deflection $w_{3}$ and the $2\times2$-matrix second-order of
differential operators%
\begin{equation}
\mathcal{L}^{\prime}\left(  \nabla_{y}\right)  =\mathcal{D}^{\prime}\left(
-\nabla_{y}\right)  ^{\top}A^{0}\mathcal{D}^{\prime}\left(  \nabla_{y}\right)
\label{3.28}%
\end{equation}
for the longitudinal displacement vector $w^{\prime}=\left(  w_{1}%
,w_{2}\right)  ^{\top}.$

In the isotropic case (see formula (1.10) in \cite{ButCaNa1}) we have%
\[
A^{0}=\left(
\begin{array}
[c]{ccc}%
\lambda^{\prime}+2\mu & \lambda^{\prime} & 0\\
\lambda^{\prime} & \lambda^{\prime}+2\mu & 0\\
0 & 0 & 2\mu
\end{array}
\right)  ,\ \ \ \lambda^{\prime}=\frac{2\lambda\mu}{\lambda+2\mu}%
\]
while $\mathcal{L}_{3}\left(  \nabla_{y}\right)  $ and $\mathcal{L}^{\prime
}\left(  \nabla_{y}\right)  $ are the bi-harmonic operator and the plane
Lam\'{e} operator
\[
\mathcal{L}_{3}^{\prime}\left(  \nabla_{y}\right)  =\frac{\mu}{3}\frac
{\lambda+\mu}{\lambda+2\mu}\Delta_{y}^{2},\ \ \ \ \mathcal{L}^{\prime}\left(
\nabla_{y}\right)  =-\mu\Delta_{y}\mathbb{I}_{2}-\left(  \lambda^{\prime}%
+\mu\right)  \nabla_{y}\nabla_{y}^{\top}.
\]

\begin{remark}
\label{remUW}Throughout the paper it is convenient to write the constructed
solutions of problems (\ref{3.7}) and (\ref{3.99}) with $f^{0}$ in the form
\begin{equation}
U^{p}\left(  y,\zeta\right)  =W^{p}\left(  \zeta,\nabla_{y}\right)  w\left(
y\right)  ,\ \ \ p=0,1,2,3, \label{3.52}%
\end{equation}
where, according to (\ref{3.8}), (\ref{3.11}) and (\ref{3.2678}),%
\begin{align}
W^{0}\left(  \zeta,\nabla_{y}\right)   &  =\left(
\begin{array}
[c]{ccc}%
0 & 0 & 0\\
0 & 0 & 0\\
0 & 0 & 1
\end{array}
\right)  ,\ \ \ W^{1}\left(  \zeta,\nabla_{y}\right)  =\left(
\begin{array}
[c]{ccc}%
1 & 0 & -\zeta\partial_{1}\\
0 & 1 & -\zeta\partial_{2}\\
0 & 0 & 0
\end{array}
\right)  ,\ \ \ \partial_{i}=\frac{\partial}{\partial y_{i}},\label{3.53}\\
W^{2}\left(  \zeta,\nabla_{y}\right)   &  =\mathbb{J}^{-1}A_{\left(
zz\right)  }^{-1}A_{\left(  zy\right)  }\left(  -\zeta\mathbb{I}_{3}%
,2^{1/2}\left(  \frac{\zeta^{2}}{2}-\frac{1}{24}\right)  \mathbb{I}%
_{3}\right)  \mathcal{D}\left(  \nabla_{y}\right)  .\nonumber
\end{align}
A concrete formula for $W^{3}\left(  \zeta,\nabla_{y}\right)  $ will not be
applied. We emphasize that (\ref{3.52}) is a linear combination of derivatives
of $w^{\prime}=\left(  w_{1},w_{2}\right)  $ of order $p-1$ and $w_{3}$ of
order $p.$ $\blacksquare$
\end{remark}

\begin{proposition}
\label{prop3.2}Let%
\begin{equation}
g_{3}=g_{3}^{0}-\sum\limits_{i=1}^{2}\frac{\partial}{\partial y_{i}}g_{3}%
^{i},\ \ \ g_{3}^{0},\ g_{3}^{i},\ g_{i}\in L^{2}\left(  \omega\right)  .
\label{3.31}%
\end{equation}
Problem (\ref{3.15}), (\ref{3.23}), (\ref{3.24}) admits a unique solution in
the energy space $\mathcal{H}=H_{0}^{1}\left(  \omega;\partial\omega\right)
^{2}\times H_{0}^{2}\left(  \omega;\partial\omega\cup\mathcal{O}\right)  .$
Moreover, this solution falls into $H^{2}\left(  \omega\right)  ^{2}\times
H^{3}\left(  \omega\right)  $ and meets estimate%
\begin{equation}
||w;H^{2}\left(  \omega\right)  ^{2}\times H^{3}\left(  \omega\right)  ||\leq
c\mathcal{N} \label{3.32}%
\end{equation}
where $\mathcal{N}$\ is the sum of norms of functions indicated in
(\ref{3.31}). The space $\mathcal{H}$\ consists of vector functions
(\ref{3.9}) in the direct product $H^{1}\left(  \omega\right)  ^{2}\times
H^{2}\left(  \omega\right)  $ of Sobolev spaces such that the Dirichlet
(\ref{3.23}) and Sobolev (\ref{3.24}) conditions are satisfied.
\end{proposition}

\textbf{Proof.} The variational formulation of problem (\ref{3.15}),
(\ref{3.23}), (\ref{3.24}) reads: to find $w\in\mathcal{H}$ fulfilling the
integral identity
\begin{equation}
\big(\mathcal{AD}(\nabla_{y})w,\mathcal{D}(\nabla_{y})v\big)_{\omega}=\left(
g_{3}^{0},v_{3}\right)  _{\omega}+\sum\limits_{i=1}^{2}\left(  \left(
g_{i},v_{i}\right)  _{\omega}+\Big(g_{3}^{i},\frac{\partial v_{3}}{\partial
y_{i}}\Big)_{\omega}\right)  ,\qquad\forall v\in\mathcal{H}. \label{3.33}%
\end{equation}
The latter, as usual, is obtained by multiplying system (\ref{3.15}) scalarly
with a test function $v\in\mathcal{H}$ and integrating by parts in $\omega$
with the help of conditions (\ref{3.23}), (\ref{3.24}) for $v$ and the formula
for $g_{3}$ in (\ref{3.31}). Clearly, the right-hand side of (\ref{3.33}) is a
continuous linear functional in $\mathcal{H}\ni v.$ The left-hand side of
(\ref{3.33}) may be chosen as a scalar product in $\mathcal{H}.$ Indeed, the
necessary properties of this bi-linear form are inherited from the positive
definiteness of the matrix $\mathcal{A}^{0}$ together with the following
Friedrichs and Korn inequalities:%
\begin{align}
||u;L^{2}\left(  \omega\right)  ||  &  \leq c_{\omega}||\nabla_{y}%
u;L^{2}\left(  \omega\right)  ||,\ \ \ \forall u\in H_{0}^{1}\left(
\omega;\partial\omega\right)  ,\label{3.34}\\
||\nabla_{y}w^{\prime};L^{2}\left(  \omega\right)  ||  &  \leq c_{\omega
}||\mathcal{D}^{\prime}\left(  \nabla_{y}\right)  w^{\prime};L^{2}\left(
\omega\right)  ||,\ \ \ \forall w^{\prime}\in H_{0}^{1}\left(  \omega
;\partial\omega\right)  ^{2}.\nonumber
\end{align}
Equivalently, we may write equation \eqref{3.33} as the minimization problem
\[
\min\left\{  \frac12\big(\mathcal{AD}(\nabla_{y})w,\mathcal{D}(\nabla
_{y})w\big)_{\omega}-(g,w)_{\omega}\ :\ w\in\mathcal{H}\right\}
\]
where $g=(g_{1},g_{2},g_{3})$ is given by \eqref{3.31}. Notice that the
Sobolev embedding theorem $H^{2}\subset C$ in $\mathbb{R}^{2}$ assures that
$\mathcal{H}$ is a closed subspace in $H^{1}\left(  \omega\right)  ^{2}\times
H^{2}\left(  \omega\right)  $ and the second inequality in (\ref{3.34}) can be
easily derived by integration by parts because $w^{\prime}$ vanishes at
$\partial\omega$.

Now the Riesz representation theorem guarantees the existence of a unique weak
solution to problem (\ref{3.33}) and the estimate%
\[
||w;H^{1}\left(  \omega\right)  ^{2}\times H^{2}\left(  \omega\right)  ||\leq
c\mathcal{N}.
\]
Finally, a result in \cite[Ch.2]{LiMa} on lifting smoothness of solutions to
elliptic problems, see also \cite{BuNa} for details about the fourth-order
equation, concludes with the inclusion $w\in H^{2}\left(  \omega\right)
^{2}\times H^{3}\left(  \omega\right)  $ and inequality (\ref{3.32}).
$\blacksquare$

\bigskip

A classical assertion, cf. \cite{Ciarlet1, Ciarlet2, LeDret, SaPa} and others,
was outlined in Theorem 7 of \cite{ButCaNa1}, namely the rescaled
displacements $h^{3/2}u_{3}\left(  h,y,h\zeta\right)  $ and $h^{1/2}%
u_{i}\left(  h,y,h\zeta\right)  ,\ i=1,2,$ taken from the three-dimensional
problem (\ref{B1.11})-(\ref{B1.14}), converge in $L^{2}\left(  \omega
\times\left(  -1/2,1/2\right)  \right)  $ as $h\rightarrow+0$ to the functions
$w_{3}\left(  y\right)  $\ and\ $w_{i}\left(  y\right)  -\zeta\frac{\partial
w_{3}}{\partial y_{i}}\left(  y\right)  ,\ $respectively, where $w=\left(
w_{1},w_{2},w_{3}\right)  ^{\top}$ is the solution of the two-dimensional
Dirichlet-Sobolev problem given in Proposition \ref{prop3.2}. Moreover, the
convergence rate $O(|\ln h|^{-1/2}),$ see Theorem \ref{th4.52} and Remark
\ref{remBoLa}, is rather slow and in the next section we will improve the
asymptotic result by constructing three-dimensional boundary layers studied in
\cite{ButCaNa1, na135, na243}.

\section{Constructing and justifying the asymptotics\label{sect4}}

\subsection{Matching the outer and inner expansions\label{sect4.1}}

Intending to apply the method of matched asymptotic expansions, cf.
\cite{Ilin, VanDyke} and \cite[Ch.2]{MaNaPl}, we regard (\ref{3.3}) as the
\textit{outer} expansion suitable at a distance from the small support area
(\ref{B1.3}). In the vicinity of the set $\theta_{h}$ we, as mentioned in
Section 3 of \cite{ButCaNa1}, use the stretched coordinates
\begin{equation}
\xi=(\eta,\zeta)=(h^{-1}y,h^{-1}z) \label{B3.35}%
\end{equation}
and engage the \textit{inner} expansion in the form%
\begin{equation}
u\left(  h,x\right)  =h^{1/2}\mathcal{P}\left(  \xi\right)  a+... \label{4.1}%
\end{equation}
where $\mathcal{P}$\ is the elastic logarithmic potential, see Section 3.4 in
\cite{ButCaNa1}, and $a\in\mathbb{R}^{4}$ is a column to be determined. The
decompositions for this potential derived in \cite{ButCaNa1} show that, when
$\rho\rightarrow+\infty,$ we have%
\begin{align}
h^{1/2}\mathcal{P}\left(  \xi\right)  a  &  =h^{1/2}\Big(\sum\limits_{p=0}%
^{3}h^{p}W^{p}\left(  \zeta,\nabla_{\eta}\right)  \Phi^{\sharp}\left(
\eta\right)  +d^{\sharp}\left(  \eta,\zeta\right)  C^{\sharp}%
\Big)a+...=\label{4.2}\\
&  =\sum\limits_{p=0}^{3}h^{p-3/2}W^{p}\left(  \zeta,\nabla_{y}\right)
\left(  \Phi^{\sharp}\left(  y\right)  -\Psi\ln h+d^{\sharp}\left(
\eta,0\right)  C^{\sharp}\right)  a+...\nonumber
\end{align}
Here, the operators $W^{p}(\zeta,\nabla_{\eta})$ were introduced in Remark
\ref{remUW} and dots stand for lower-order asymptotic terms and we have taken
into account formulas%
\begin{align}
W^{p}\left(  \zeta,\nabla_{\eta}\right)   &  =h^{p-1}W^{p}\left(  \zeta
,\nabla_{y}\right)  H,\ \ \ H=\mathrm{diag}\{1,1,h\},\label{4.3}\\
d^{\sharp}\left(  \eta,\zeta\right)   &  =H^{-1}d^{\sharp}\left(
y,\zeta\right)  ,\ \Phi^{\sharp}\left(  \eta\right)  =H^{-1}\left(
\Phi^{\sharp}\left(  y\right)  -d^{\sharp}\left(  y,0\right)  \Psi\ln h\right)
\nonumber
\end{align}
where $d^{\sharp}$ and $\Phi^{\sharp}$ are the following $3\times4$-matrices
\begin{align}
d^{\sharp}\left(  \eta,\zeta\right)   &  =\left(
\begin{array}
[c]{cccc}%
1 & 0 & 0 & \zeta\\
0 & 1 & -\zeta & 0\\
0 & 0 & \eta_{2} & -\eta_{1}%
\end{array}
\right) \label{B3.49}\\
\Phi^{\sharp}\left(  y\right)   &  =\left(  d^{\sharp}\left(  -\nabla
_{y},0\right)  ^{\top}\Phi\left(  y\right)  ^{\top}\right)  ^{\top}=\left(
\begin{array}
[c]{c}%
\Phi^{\prime}\left(  y\right) \\
\\
0\quad\quad0
\end{array}%
\begin{array}
[c]{cc}%
0 & 0\\
0 & 0\\
\Phi_{3}^{2}\left(  y\right)  & \Phi_{3}^{1}\left(  y\right)
\end{array}
\right)  ,\ \ \Phi_{3}^{j}=\frac{\partial\Phi_{3}}{\partial y_{j}},\nonumber
\end{align}
which are composed, respectively, from rigid motions and the nondegenerate
block-diagonal $4\times4$-matrix%
\begin{equation}
\Psi=\operatorname*{diag}\left\{  \Psi^{\prime},\Psi_{3},\Psi_{3}\right\}
\label{4.311}%
\end{equation}
where $\Psi^{\prime}$ is a non-degenerate numeral $2\times2$-matrix and
$\Psi_{3}$ is a scalar in the representations%
\begin{equation}
\Phi^{\prime}(y)=\Psi^{\prime}\ln r+\psi^{\prime}(\varphi),\ \ \ \Phi
_{3}(y)=r^{2}\left(  -\frac{1}{2}\Psi_{3}\ln r+\psi_{3}(\varphi)\right)
\label{B346345}%
\end{equation}
of the fundamental matrix of the differential operator (\ref{3.28}) of size
$2\times2$ and the fundamental solution of the scalar differential operator
(\ref{3.27}). These were described in detail in Section 3.4 of \cite{ButCaNa1}
on the base of general results in \cite{GelShilov}. Moreover, $(r,\varphi)$ is
the polar coordinate system and $\psi^{\prime},$ $\psi_{3}$ are smooth on the
unit circle $\mathbb{S}^{1}\ni\varphi.$ Furthermore, we have used the diagonal
matrix $H=\operatorname*{diag}\left\{  1,1,h\right\}  $ and the obvious
formulas%
\begin{equation}
d^{\sharp}\left(  \eta,\zeta\right)  =\left(  W^{0}+W^{1}\left(  \zeta
,\nabla_{\eta}\right)  \right)  d^{\sharp}\left(  \eta,0\right)
,\ W^{2}\left(  \zeta,\nabla_{\eta}\right)  d^{\sharp}\left(  \eta,0\right)
=W^{3}\left(  \zeta,\nabla_{\eta}\right)  d^{\sharp}\left(  \eta,0\right)  =0.
\label{4.4}%
\end{equation}
Note that $\ln h$ comes to (\ref{4.3}) from the relation $\ln\left\vert
\eta\right\vert =\ln\left\vert y\right\vert -\ln h$ and $H$ is caused by
different orders of differentiation in $W^{p}\left(  \zeta,\nabla_{\eta
}\right)  $ and different degrees of $\rho=\left\vert \eta\right\vert $ in
matrix functions but $H$ disappears from the final expression in (\ref{4.2}).

We also will need the following expansion of the elastic logarithmic
potential, see Section 3.6 in \cite{ButCaNa1},%
\begin{equation}
\mathcal{P}(\xi)=\left(  1-\chi_{\theta}(\eta)\right)  \Bigg(\sum
\limits_{p=0}^{3}W^{p}\left(  \zeta,\nabla_{\eta}\right)  \Phi^{\sharp}\left(
\eta\right)  +d^{\sharp}\left(  \eta,0\right)  C^{\sharp}+\Upsilon^{\sharp
}\left(  \eta\right)  \Bigg)+\widetilde{\widetilde{\mathcal{P}}}(\xi)
\label{8.8N}%
\end{equation}
which had been used in (\ref{4.2}). Ingredients of (\ref{8.8N}) were defined
in (\ref{3.53}) and (\ref{B3.49})-(\ref{B346345}) while $C^{\sharp}$ stands
for the elastic capacity matrix, $\chi_{\theta}$ a smooth compactly supported
function which equals 1 in the vicinity of the clamped area $\theta,$ the
vector function $\Upsilon^{\sharp}$ was described in Section 3.6 of
\cite{ButCaNa1} but it plays an auxiliary role and an explicit form is not
used below, and the remainders $\widetilde{\widetilde{\mathcal{P}}}$ admits
the estimates for a big $\rho=|\eta|$%
\begin{equation}
|\nabla_{\eta}^{p}\partial_{\zeta}^{q}\widetilde{\widetilde{\mathcal{P}}}%
(\xi)|\leq c_{pq}\rho^{-1+\varepsilon},\ |\nabla_{\eta}^{p}\partial_{\zeta
}^{q}\widetilde{\widetilde{\mathcal{P}}}_{3}(\xi)|\leq c_{pq}\rho
^{\varepsilon},\ p,q=0,1,2,...,\ \varepsilon>0. \label{8.9N}%
\end{equation}
The presence of $\Phi^{\sharp}\left(  y\right)  $ drives all detached terms in
(\ref{4.3}) from the space $H^{1}\left(  \omega\right)  ^{2}\times
H^{2}\left(  \omega\right)  $ where the solution $w$ of the limit problem
(\ref{3.15}), (\ref{3.23}), (\ref{3.24}) was found in Proposition
\ref{prop3.2}. On the other hand, the Neumann boundary condition (\ref{B1.12})
on the lower base $\Sigma^{-}$ does not hold on the small set $\overline
{\theta_{h}}=\Sigma^{-}\setminus\Sigma^{\bullet}$ which shrinks to the
coordinate origin $\mathcal{O}$. Hence, the dimension reduction procedure in
Section \ref{sect3} was not able to deduce the differential equations at the
point $\mathcal{O}$\ and we did not have a regular reason to cast out
solutions with reasonable singularities. Let us introduce such solutions.

We denote by $G^{\prime}$ the Green $2\times2$-matrix of the Dirichlet problem
in $\omega$ for the matrix operator (\ref{3.28}). Its particular value
$G^{\prime}\left(  \cdot,\mathcal{O}\right)  $ is a distributional solution of
the problem%
\begin{equation}
\mathcal{L}^{\prime}\left(  \nabla_{y}\right)  G^{\prime}\left(
y,\mathcal{O}\right)  =\delta\left(  y\right)  \mathbb{I}_{2},\ y\in
\omega,\ \ \ \ G^{\prime}\left(  y,\mathcal{O}\right)  =0,\ y\in\partial
\omega, \label{4.5}%
\end{equation}
where $\delta$ is the Dirac mass and $\mathbb{I}_{2}$\ is the unit $2\times2$-matrix.

As known, the Green function $G_{3}$ of the Dirichlet problem for the
fourth-order scalar differential operator (\ref{3.27}) belongs to
$H^{2}\left(  \omega\right)  \subset C\left(  \omega\right)  $ and its value%
\[
G_{3}\left(  \mathcal{O},\mathcal{O}\right)  =\left(  \mathcal{A}%
_{3}\mathcal{D}_{3}\left(  \nabla_{y}\right)  G_{3}\left(  \cdot
,\mathcal{O}\right)  ,\mathcal{D}_{3}\left(  \nabla_{y}\right)  G_{3}\left(
\cdot,\mathcal{O}\right)  \right)  _{\omega}%
\]
is positive, see, e.g., \cite{BuNa}. The derivative in the second argument
\begin{equation}
G_{3,i}\left(  y,\mathcal{O}\right)  =\frac{\partial}{\partial\mathbf{y}_{i}%
}G_{3}\left(  y,\mathbf{y}\right)  |_{\mathbf{y}=0} \label{4.00}%
\end{equation}
leaves the space $H^{2}\left(  \omega\right)  $ but remains
H\"{o}lder-continuous in $\omega$ and satisfies the equation%
\begin{equation}
\mathcal{L}_{3}\left(  \nabla_{y}\right)  G_{3.i}\left(  y,\mathcal{O}\right)
=-\frac{\partial}{\partial y_{i}}\delta\left(  y\right)  ,\ y\in\omega,
\label{4.6}%
\end{equation}
and the Dirichlet conditions (\ref{3.23}). However, the Sobolev condition
(\ref{3.24}) is not achieved yet and we set%
\begin{equation}
G_{3}^{i}\left(  y,\mathcal{O}\right)  =G_{3,i}\left(  y,\mathcal{O}\right)
-G_{3,i}\left(  \mathcal{O},\mathcal{O}\right)  G_{3}\left(  \mathcal{O}%
,\mathcal{O}\right)  ^{-1}G_{3}\left(  y,\mathcal{O}\right)  , \label{4.7}%
\end{equation}
cf. \cite{BuNa}, so that $G_{3}^{i}\left(  \mathcal{O},\mathcal{O}\right)
=0$\ and\ $G_{3}^{i}\left(  y,\mathcal{O}\right)  =\partial_{n}G_{3}%
^{i}\left(  y,\mathcal{O}\right)  =0,\ y\in\partial\omega.$ Finally we
introduce the matrix function of size $3\times4$%
\[
G^{\sharp}\left(  y\right)  =\left(
\begin{array}
[c]{c}%
G^{\prime}\left(  y,\mathcal{O}\right) \\
\\
0\quad\quad0
\end{array}%
\begin{array}
[c]{cc}%
0 & 0\\
0 & 0\\
-G_{3}^{2}\left(  y,\mathcal{O}\right)  & G_{3}^{1}\left(  y,\mathcal{O}%
\right)
\end{array}
\right)
\]
which, by virtue of (\ref{4.5}) and (\ref{4.6}), admits the representation%
\begin{equation}
G^{\sharp}\left(  y\right)  =\Phi^{\sharp}\left(  y\right)  +\widehat
{G}^{\sharp}\left(  y\right)  \label{4.8}%
\end{equation}
where $\Phi^{\sharp}$ is the singular matrix function in (\ref{B3.49}) and
$\widehat{G}^{\sharp}$ is the regular part. According to general properties of
the Green matrices of second-order systems, see \cite{GelShilov}, first two
lines $\widehat{G}_{i}^{\sharp},$ $i=1,2,$ in $\widehat{G}^{\sharp}$ involve
smooth functions in $\overline{\omega},$ however the singularity $O\left(
r^{2}\left\vert \ln r\right\vert \right)  $ of the subtrahend in (\ref{4.7})
moves some entries of $\widehat{G}_{3}^{\sharp}$ from $C^{2}\left(
\omega\right)  $ but keep them still in the H\"{o}lder space $C^{1,\alpha
}\left(  \omega\right)  $ with any $\alpha\in\left(  0,1\right)  .$ Hence, the
remainder in (\ref{4.8}) verifies the relations%
\begin{equation}
\widehat{G}^{\sharp}\left(  y\right)  =d^{\sharp}\left(  y,0\right)
\mathcal{G}^{\sharp}+\widetilde{G}^{\sharp}\left(  y\right)
,\ \ \ \ \ \widetilde{G}_{i}^{\sharp}\left(  y\right)  =O\left(  r\right)
,\ \ i=1,2,\ \ \ \widetilde{G}_{3}^{\sharp}\left(  y\right)  =O(r^{2}(1+|\ln
r|)), \label{4.9}%
\end{equation}
where $\mathcal{G}^{\sharp}=\mathcal{G}^{\sharp}\left(  A,\omega\right)  $ is
a numeric matrix of size $4\times4.$ Properties of this matrix are just the
same as ones of the elastic logarithmic capacity and verification of its
symmetry is a sufficient simplification of the proof of Theorem 13 in
\cite{ButCaNa1}.

\begin{lemma}
\label{lem4.1}The $4\times4$-matrix $\mathcal{G}^{\sharp}=\mathcal{G}^{\sharp
}\left(  A,\omega\right)  $ in (\ref{4.9}) is symmetric.
\end{lemma}

Let us put into the outer expansion (\ref{3.3}) the singular solution%
\begin{equation}
w\left(  y\right)  =\widehat{w}(y)+G^{\sharp}(y)a \label{4.11}%
\end{equation}
of problem (\ref{3.15}), (\ref{3.23}), (\ref{3.24}) with the same coefficient
column $a\in\mathbb{R}^{4}$ as in (\ref{4.1}) and the smooth solution
$\widehat{w}\in H^{2}\left(  \omega\right)  ^{2}\times H^{3}\left(
\omega\right)  $ given in Proposition \ref{prop3.2}. With a reference to the
Sobolev embedding theorem $H^{2}\subset C$ and the standard Hardy inequalities
(see, e.g., Section 2.1 in \cite{ButCaNa1}), we write%
\begin{equation}
\widehat{w}\left(  y\right)  =d^{\sharp}(y,0)F+\widetilde{w}(y) \label{4.12}%
\end{equation}
where the coefficient column $F\in\mathbb{R}^{4}$ and the remainder
$\widetilde{w}$ satisfy the estimates
\begin{gather}
\left\vert F\right\vert \leq c\mathcal{N},\label{4.13}\\
\int_{\omega}r^{-2}(1+\left\vert \ln r\right\vert )^{-2}\Big(|\nabla
_{y}\widetilde{w}^{\prime}(y)|^{2}+r^{-2}|\widetilde{w}^{\prime}%
(y)|^{2}+|\nabla_{y}^{2}\widetilde{w}_{3}(y)|^{2}\label{4.14}\\
+r^{-2}|\nabla_{y}\widetilde{w}_{3}(y)|^{2}+r^{-4}|\widetilde{w}_{3}%
(y)|^{2}\Big)dy\leq c\mathcal{N},\nonumber
\end{gather}
while $\mathcal{N}$\ is the sum of norms of functions in (\ref{3.31}).

Notice that, according to definition of the Green matrix, we have%
\[
F=\int_{\omega}G^{\sharp}\left(  y\right)  g\left(  y\right)  dy.
\]
Using formulas (\ref{4.11}), (\ref{4.12}) and notation (\ref{3.52}), we obtain
for the outer expansion (\ref{3.3}) that%
\begin{align}
h^{-3/2}\sum\limits_{p=0}^{3}h^{p}U^{p}\left(  y,\zeta\right)   &
=h^{-3/2}\sum\limits_{p=0}^{3}h^{p}W^{p}\left(  \zeta,\nabla_{y}\right)
w\left(  y\right) \label{4.15}\\
&  =\sum\limits_{p=0}^{3}h^{p-3/2}W^{p}\left(  \zeta,\nabla_{y}\right)
\left(  d^{\sharp}\left(  y,0\right)  F+\left(  \Phi^{\sharp}\left(  y\right)
+d^{\sharp}\left(  y,0\right)  \right)  a\right)  +...\nonumber
\end{align}
where dots again stand for lower-order asymptotic terms. The method of matched
asymptotic expansions requires that the expressions detached in (\ref{4.15})
and (\ref{4.2}), coincide with each other. This coincidence implies the system
of linear algebraic equations%
\begin{equation}
F+\mathcal{G}^{\sharp}a=-\ln h~\Psi a+C^{\sharp}a \label{4.16}%
\end{equation}
and, thus,
\begin{equation}
a\left(  \ln h\right)  =\left(  \left\vert \ln h\right\vert \Psi+C^{\sharp
}\left(  A,\theta\right)  -\mathcal{G}^{\sharp}\left(  A,\omega\right)
\right)  ^{-1}F. \label{4.17}%
\end{equation}
We here display the dependence of $C^{\sharp}$ and $\mathcal{G}^{\sharp}$ on
the stiffness matrix $A$ and the domains $\theta$ and $\omega.$ Since the
matrix (\ref{4.311}) is not degenerate, the matrix%
\begin{equation}
M^{\sharp}\left(  \ln h\right)  =\left\vert \ln h\right\vert \Psi+C^{\sharp
}\left(  A,\theta\right)  -\mathcal{G}^{\sharp}\left(  A,\omega\right)
\label{4.18}%
\end{equation}
is invertible for a small $h\in(0,1)$, too. Moreover, column (\ref{4.17}) is a
rational vector function in $\left\vert \ln h\right\vert =-\ln h$ while
clearly%
\begin{equation}
a\left(  \ln h\right)  =\left\vert \ln h\right\vert ^{-1}\Psi^{-1}%
F+O(\left\vert \ln h\right\vert ^{-2}). \label{4.19}%
\end{equation}

Formula (\ref{4.17}) gives concrete expressions to all terms in (\ref{3.3}),
(\ref{4.11}) and (\ref{4.1}) so that our formal construction of main terms in
the outer and inner expansions is completed.

\subsection{The final assumptions on the right-hand sides\label{sect4.2}}

To provide the necessary properties of the regular solution (\ref{4.12}) of
problem (\ref{3.15}), (\ref{3.23}), (\ref{3.24}), we put the following
requirement on terms in representation (\ref{3.1}):%
\begin{equation}
f^{0}\in L^{2}\left(  \omega\times\left(  -\tfrac{1}{2},\tfrac{1}{2}\right)
\right)  ^{3},\ \ \ f_{3}^{1}\in H^{-1}\left(  \omega\right)  . \label{FF1}%
\end{equation}
The latter inclusion, as usual, cf. \cite{LiMa}, means that
\begin{equation}
f_{3}^{1}\left(  y\right)  =\nabla_{y}^{\top}\overrightarrow{f}_{3}^{1}\left(
y\right)  +f_{30}^{1}\left(  y\right)  ,\ \ \ \ \overrightarrow{f}_{3}%
^{1}=\left(  f_{31}^{1},f_{32}^{1}\right)  ^{\top}\in L^{2}\left(
\omega\right)  ^{2},\ \ \ \ f_{30}^{1}\in L^{2}\left(  \omega\right)  .
\label{FF2}%
\end{equation}
More precisely, the right-hand side (\ref{3.22}) of the third equation in
system (\ref{3.15}) gives rise to the continuous functional%
\begin{equation}
\left(  g_{3},\mathsf{v}_{3}\right)  _{\omega}:=\left(  f_{30}^{1}%
,\mathsf{v}_{3}\right)  _{\omega}-\sum\limits_{i=1}^{2}\left(  f_{3i}^{1}%
+\int_{-1/2}^{1/2}\zeta f_{i}^{0}\left(  \cdot,\zeta\right)  d\zeta
,\frac{\partial\mathsf{v}_{3}}{\partial y_{i}}\right)  _{\omega}%
,\ \ \ \mathsf{v}_{3}\in H_{0}^{1}\left(  \omega\right)  . \label{FF3}%
\end{equation}
In view of (\ref{3.17}), (\ref{FF1}) and (\ref{FF3}), (\ref{FF2}) condition
(\ref{3.31}) in Proposition \ref{prop3.2} is fulfilled so that the solution
$\widehat{w}\in H^{2}\left(  \omega\right)  ^{2}\times H^{3}\left(
\omega\right)  $ meets estimate (\ref{3.32}) where, as well as in
(\ref{4.13}), (\ref{4.14}), $\mathcal{N}$ can be now fixed as%
\begin{equation}
\mathcal{N}=\left\Vert f^{0};L^{2}\left(  \omega\times\left(  -\tfrac{1}%
{2},\tfrac{1}{2}\right)  \right)  \right\Vert +\sum\limits_{p=0}^{2}\left\Vert
f_{3p}^{1};L^{2}\left(  \omega\right)  \right\Vert . \label{FF4}%
\end{equation}
For the remainder $\widetilde{f}$ in (\ref{3.1}), we suppose that the
expression%
\begin{equation}
\widetilde{\mathcal{N}}=h^{-1/2}|\ln h|^{-1}\Big(h^{-1}||s_{h}^{2}S_{h2}%
^{-1}\widetilde{f}_{3};L^{2}\left(  \Omega_{h}\right)  ||+\sum\limits_{i=1}%
^{2}||s_{h}S_{h1}^{-1}\widetilde{f}_{i};L^{2}\left(  \Omega_{h}\right)
||\Big) \label{FF5}%
\end{equation}
gets the same order as the expression (\ref{FF4}). Notice that weights
\begin{equation}
s_{h}(y)=h+\mathrm{dist}(y,\partial\omega),\ S_{hq}(y)=(h^{2}+|y|^{2}%
)^{-q/2}(1+|\ln(h^{2}+\left\vert y\right\vert ^{2})|)^{-1},\ q=0,1,
\label{B2425}%
\end{equation}
come from the anisotropic Korn inequality (2.10) in \cite{ButCaNa1} and make
the norms on the right-hand side of (\ref{FF5}) much weaker that the common
Lebesgue norms. The factor $h^{-1/2}\left\vert \ln h\right\vert ^{-1}$ is
consistent with the final error estimate in Theorem \ref{th4.51}. In other
words, our way to signify the smallness of this remainder does not spoil the
precision of the two-dimensional model.

\subsection{The global approximation solution\label{sect4.3}}

We proceed with constructing a vector function $\mathbf{u}\in H_{0}^{1}\left(
\Omega_{h};\upsilon_{h}\cup\theta_{h}\right)  ^{3}$ which is a bit cumbersome
but very convenient for an estimation of the difference $u-\mathbf{u}$ between
the true and approximate solutions of problem (\ref{B1.11})-(\ref{B1.14}). In
the next section we will get rid of unnecessary terms in $\mathbf{u}$\ to
conclude a fair asymptotics of $u.$

Although in the previous section we have used the method of matched
expansions, we now apply an asymptotic structure attributed to the method of
compound expansions, see a comparison of the methods in \cite[Ch.2]{MaNaPl}.
We give priority to expansion (\ref{4.1}) and insert into (\ref{3.3}) the
decaying near $\mathcal{O}$ component $\widetilde{w},$ see (\ref{4.12}) and
(\ref{4.14}), instead of the whole singular solution (\ref{4.11}). We also
introduce the cut-off functions%
\begin{gather}
X_{h}^{\omega}\left(  y\right)  =1-\chi\left(  h^{-1}\left\vert n\right\vert
\right)  ,\ \ \ X_{h}^{\theta}\left(  y\right)  =1-\chi\left(  2h^{-1}%
r/R_{\theta}\right) \label{4.20}\\
\chi\in C^{\infty}\left(  R\right)  ,\ \chi(r)=1\text{ \ for }r<\tfrac{1}%
{2},\ \ \chi(r)=0\text{ \ for }r\geq1,\ 0\leq\chi\leq1,\nonumber
\end{gather}
which help to fulfill the Dirichlet conditions (\ref{B1.14}) and (\ref{B1.13})
because $X_{h}^{\omega}$ and $X_{h}^{\theta}$ vanish near the clamped sets
$\upsilon_{h}$ and $\theta_{h}$ respectively. To this end, $R_{\theta}$ is
fixed such that $\overline{\theta}$ belongs to the disk $\mathbb{B}_{R}^{2}$
of radius $R=R_{\theta}$.

We set%
\begin{align}
\mathbf{u}\left(  h,x\right)   &  =h^{-1/2}X_{h}^{\omega}\left(  y\right)
P\left(  h^{-1}x\right)  a\left(  \ln h\right)  +\label{4.21}\\
&  +h^{-3/2}X_{h}^{\omega}\left(  y\right)  \sum\limits_{p=0}^{2}h^{p}%
W^{p}\left(  \zeta,\nabla_{y}\right)  \left(  X_{h}^{\theta}\left(  y\right)
\widetilde{w}\left(  y\right)  -h\mathbf{w}\left(  h,y\right)  \right)
\nonumber
\end{align}
where the following correction term is introduced%
\begin{equation}
\mathbf{w}\left(  h,y\right)  =\left(  1-\chi\left(  R_{\omega}^{-1}r\right)
\right)  \Upsilon^{\sharp}\left(  h^{-1}y\right)  \label{4.22}%
\end{equation}
with $\Upsilon^{\sharp}$ taken from (\ref{8.8N}) and $R_{\omega}>0$ is such
that $\overline{\mathbb{B}_{R_{\omega}}}\subset\omega.$ We emphasize that,
according to definition of the differential operators (\ref{3.53}) in
(\ref{3.52}), the inclusion $\widetilde{w}\in H^{2}\left(  \omega\right)
^{2}\times H^{3}\left(  \omega\right)  $ assures that%
\begin{equation}
W^{p}\left(  \zeta,\nabla_{y}\right)  \widetilde{w}\in H^{3-p}\left(
\omega\right)  ^{3}\text{ for any }\zeta\in\left(  -1/2,1/2\right)
\label{4.000}%
\end{equation}
and, therefore, the term $W^{3}\left(  \zeta,\nabla_{y}\right)  \widetilde
{w}\left(  y\right)  $ does not belong to the energy space $H^{1}\left(
\Omega_{h}\right)  ^{3}$ and is excluded from the global approximation
solution, i.e., $p=0,1,2$ but $p\neq3$ in (\ref{4.21}).

Thanks to the cut-off functions (\ref{4.20}) and the Dirichlet conditions for
$\mathcal{P}$, the displacement field satisfies conditions (\ref{B1.13}),
(\ref{B1.14}). Hence, the difference $\mathbf{v}=u-\mathbf{u\in}\mathring
{H}^{1}\left(  \Omega_{h};\upsilon_{h}\cup\theta_{h}\right)  ^{3}$ can be
taken as a test function in the integral identity (\ref{B1.16}). Subtracting
from both sides the scalar product $\left(  AD\left(  \nabla\right)
\mathbf{u},D\left(  \nabla\right)  \mathbf{v}\right)  _{\Omega_{h}},$ we
arrive at the formula%
\begin{equation}
\left(  AD\left(  \nabla\right)  \mathbf{v},D\left(  \nabla\right)
\mathbf{v}\right)  _{\Omega_{h}}=\left(  f,\mathbf{v}\right)  _{\Omega_{h}%
}-\left(  AD\left(  \nabla\right)  \mathbf{u},\mathbf{v}\right)  _{\Omega_{h}%
}. \label{4.24}%
\end{equation}
By the Korn inequality (see (2.10) and Theorem 2 in \cite{ButCaNa1}), the
left-hand side of (\ref{4.24}) exceed the product $c_{A}|||v;\Omega
_{h}|||_{\bullet}^{2}$ with the weighted anisotropic Sobolev norm
\begin{gather}
|||u;\Omega_{h}|||_{\bullet}^{2}=\int_{\Omega_{h}}\Bigg(\sum_{i=1}%
^{2}\Big(\left\vert \nabla_{y}u_{i}\right\vert ^{2}+\frac{h^{2}}{s_{h}^{2}%
}S_{h1}^{2}\Big(\left\vert \frac{\partial u_{i}}{\partial z}\right\vert
^{2}+\left\vert \frac{\partial u_{3}}{\partial y_{i}}\right\vert
^{2}\Big)+\frac{1}{s_{h}^{2}}S_{h1}^{2}\left\vert u_{i}\right\vert
^{2}\Big)\label{B2.7}\\
+\left\vert \frac{\partial u_{3}}{\partial z}\right\vert ^{2}+\frac{h^{2}%
}{s_{h}^{4}}S_{h2}^{2}\left\vert u_{3}\right\vert ^{2}\Bigg)dx.\nonumber
\end{gather}
and our immediate objective becomes to examine the right-hand side. The
complication arises from the cut-off functions (\ref{4.20}) which are included
into (\ref{4.21}) and meet the relations
\begin{align}
|\nabla_{y}X_{h}^{\omega}\left(  y\right)  |  &  \leq ch^{-1}%
,\ \ \ \operatorname*{supp}|\nabla_{y}X_{h}^{\omega}|\subset\omega\left(
h\right)  :=\left\{  y\in\omega:0>n>-h\right\}  ,\label{4.25}\\
|\nabla_{y}X_{h}^{\theta}\left(  y\right)  |  &  \leq ch^{-1}%
,\ \ \ \operatorname*{supp}|\nabla_{y}X_{h}^{\theta}|\subset\overline
{\mathbb{B}_{2R_{\omega}h}}.\nonumber
\end{align}
Denoting $\mathbf{v}^{\omega}=X_{h}^{\omega}\mathbf{v}$ and $\widetilde
{w}^{\theta}=X_{h}^{\theta}\widetilde{w},$ we obtain%
\begin{align*}
||D\left(  \nabla\right)  \mathbf{v}^{\omega};L^{2}\left(  \Omega_{h}\right)
||  &  \leq c\left(  ||D\left(  \nabla\right)  \mathbf{v};L^{2}\left(
\Omega_{h}\right)  ||+h^{-1}||\mathbf{v};L^{2}\left(  \omega\left(  h\right)
\times\left(  -h/2,h/2\right)  \right)  ||\right) \\
&  \leq C||D\left(  \nabla\right)  \mathbf{v};L^{2}\left(  \Omega_{h}\right)
||=:C\mathbf{n},\\
||\widetilde{w}^{\theta};H^{2}\left(  \omega\right)  ^{2}\times H^{3}\left(
\omega\right)  ||  &  \leq c\Big(||\widetilde{w};H^{2}\left(  \omega\right)
^{2}\times H^{3}\left(  \omega\right)  ||^{2}+J\Big)\\
&  \leq c||\widetilde{w};H^{2}\left(  \omega\right)  ^{2}\times H^{3}\left(
\omega\right)  ||^{2}\leq C\mathcal{N},
\end{align*}
where $J$ is the weighted integral on the left in (\ref{4.14}). In both
inequalities we have taken into account that the weights $s_{h}\left(
y\right)  $ and $r\left(  1+\left\vert \ln r\right\vert \right)  $ are small
in the boundary strip $\omega\left(  h\right)  $ and the disk $\mathbb{B}%
_{2R_{\theta}h}.$

\subsection{Estimating the discrepancies\label{sect4.4}}

Let $h^{-3/2}X_{h}^{\omega}\mathcal{W}-h^{-1/2}X_{h}^{\omega}\mathbf{W}$ stand
for the last sum in (\ref{4.21}). We have%
\begin{align}
\left(  AD\left(  \nabla\right)  \mathbf{u},D\left(  \nabla\right)
\mathbf{v}\right)  _{\Omega_{h}}  &  =h^{-1/2}\left(  AD\left(  \nabla\right)
\mathcal{P}a,D\left(  \nabla\right)  \mathbf{v}^{\omega}\right)  _{\Omega_{h}%
}\label{4.27}\\
&  +h^{-1/2}\left(  A\left(  D\left(  \nabla\right)  X_{h}^{\omega}\right)
\left(  \mathcal{P}a+h^{-1}\mathcal{W}-\mathbf{W}\right)  ,D\left(
\nabla\right)  \mathbf{v}\right)  _{\Omega_{h}}\nonumber\\
&  +\left(  AD\left(  \nabla\right)  \left(  \mathcal{P}a+h^{-1}%
\mathcal{W}-\mathbf{W}\right)  ,\left(  D\left(  \nabla\right)  X_{h}^{\omega
}\right)  \mathbf{v}\right)  _{\Omega_{h}}\nonumber\\
&  +h^{-1/2}\left(  AD\left(  \nabla\right)  \mathbf{W},D\left(
\nabla\right)  \mathbf{v}^{\omega}\right)  _{\Omega_{h}}+h^{-3/2}\left(
AD\left(  \nabla\right)  \mathcal{W},D\left(  \nabla\right)  \mathbf{v}%
^{\omega}\right)  _{\Omega_{h}}\nonumber\\
&  =:I_{1}+I_{2}+I_{3}+I_{4}.\nonumber
\end{align}
First of all, we observe that $\mathbf{v}=0$ at the lateral side $\upsilon
_{h}$ and, therefore, after extending $\mathbf{v}^{\omega}$ as null over the
layer $\mathbb{R}^{2}\times\left(  -h/2,h/2\right)  $ we make the coordinate
dilation (\ref{B3.35}) and integrate by parts to conclude that $I_{1}=0$ since
$\mathcal{P}$ is a solution of the homogeneous elasticity problem in the
infinite layer $\Lambda=\mathbb{R}^{2}\times\left(  -1/2,1/2\right)  $ clamped
over the set $\overline{\theta}\times\left\{  -1/2\right\}  $ on its lower base.

Next, formulas (\ref{8.8N}) and (\ref{4.11}), (\ref{4.12}) together with the
matching condition (\ref{4.16}) yield%
\begin{align*}
&  \mathcal{P}\left(  \xi\right)  a+h^{-1}\mathcal{W}\left(  h,x\right)
-\mathbf{W}\left(  h,x\right) \\
&  =\widetilde{\widetilde{\mathcal{P}}}\left(  \xi\right)  a+h^{2}\left(
W^{3}\left(  \zeta,\nabla_{y}\right)  \Phi^{\sharp}\left(  \eta\right)
+d^{\sharp}\left(  \eta,0\right)  C^{\sharp}+\Upsilon^{\sharp}\left(
\eta\right)  \right)  a+h^{-1}\sum\limits_{p=0}^{2}h^{p}W^{p}\left(
\zeta,\nabla_{y}\right)  w\left(  y\right)  .
\end{align*}
The first term on the right is estimated by means of (\ref{8.9N}) noticing
that $\rho^{-1}=hr^{-1}=O\left(  h\right)  $ in the thin boundary strip
$\omega\left(  h\right)  .$ The second term and its derivatives become
$O\left(  h^{2}\left(  1+\left\vert \ln h\right\vert \right)  \right)  $ near
the plate edge due to formulas (\ref{B3.49}), (\ref{4.4}) and (\ref{4.19}).
Furthermore,
\begin{align*}
&  \int_{\omega\left(  h\right)  }\left(  s_{h}^{-2}\left\vert w^{\prime
}\right\vert ^{2}+\left\vert \nabla_{y}w^{\prime}\right\vert ^{2}+s_{h}%
^{-4}\left\vert w_{3}\right\vert ^{2}+s_{h}^{-2}\left\vert \nabla_{y}%
w_{3}\right\vert ^{2}+\left\vert \nabla_{y}^{2}w_{3}\right\vert ^{2}\right)
dy\\
&  \leq c\int_{\omega\left(  h\right)  }\left(  \left\vert \nabla_{y}%
w^{\prime}\right\vert ^{2}+\left\vert \nabla_{y}^{2}w_{3}\right\vert
^{2}\right)  dy\\
&  \leq ch\left(  \left\Vert w^{\prime};H^{2}\left(  \omega\setminus
\mathbb{B}_{R_{\omega}/2}\right)  \right\Vert ^{2}+\left\Vert w_{3}%
;H^{3}\left(  \omega\setminus\mathbb{B}_{R_{\omega}/2}\right)  \right\Vert
^{2}\right) \\
&  \leq ch\mathcal{N}%
\end{align*}
follows from formulas (\ref{4.14}), (\ref{4.19}), the standard Hardy
inequalities (cf. Section 2.1 in \cite{ButCaNa1}) and another consequence of
the Newton-Leibnitz formula with some fixed $T\geq R_{\omega}/2,$ namely%
\begin{align*}
\int_{0}^{h}\left\vert U\left(  t\right)  \right\vert ^{2}dt  &  =\int_{0}%
^{h}\int_{t}^{T}\frac{\partial}{\partial\tau}\left(  \chi\left(
\tau/T\right)  U\left(  \tau\right)  ^{2}\right)  d\tau dt\leq c\int_{0}%
^{h}\int_{t}^{T}\left(  \left\vert \frac{dU}{d\tau}\left(  \tau\right)
\right\vert ^{2}+\left\vert U\left(  \tau\right)  \right\vert ^{2}\right)
d\tau dt\\
&  \leq ch\left\Vert U;H^{1}\left(  0,T\right)  \right\Vert ^{2}.
\end{align*}
Estimating the matrix function $D\left(  \nabla\right)  X_{h}^{\omega}$
according to (\ref{4.25}) and mentioning that, by virtue of weight
(\ref{B2425}) in norm (\ref{B2.7}),%
\[
\left\Vert \left(  D\left(  \nabla\right)  X_{h}^{\omega}\right)
\mathbf{v};L^{2}\left(  \Omega_{h}\right)  \right\Vert \leq c|||\mathbf{v}%
;\omega\left(  h\right)  \times\left(  -h/2,h/2\right)  |||_{\bullet}\leq
c\mathbf{n},
\]
we deduce that%
\begin{align}
\left\vert I_{2}\right\vert  &  \leq ch^{-1/2}h^{1/2}\big(h^{-1}h\left(
1+\left\vert \ln h\right\vert \right)  ^{2}\left\vert a\right\vert \left(
\operatorname*{mes}\nolimits_{2}\omega\left(  h\right)  \right)
^{1/2}\label{4.30}\\
&  +\big\Vert w;H^{1}\left(  \omega\left(  h\right)  \right)  ^{2}\times
H^{2}\left(  \omega\left(  h\right)  \right)  \big\Vert^{2}+h^{2}\left(
1+\left\vert \ln h\right\vert \right)  ^{2}\left\vert a\right\vert \left(
\operatorname*{mes}\nolimits_{2}\omega\left(  h\right)  \right)
^{1/2}\big)\mathbf{n}\nonumber\\
&  \leq ch^{1/2}\left\vert \ln h\right\vert \mathbf{n}.\nonumber
\end{align}
Note that, in the middle part of (\ref{4.30}), $h^{-1/2}$ is the original
common factor, $h^{1/2}$ comes from integration in $z,$ $\operatorname*{mes}%
\nolimits_{2}\omega\left(  h\right)  =O\left(  h\right)  $ and $\left(
1+\left\vert \ln h\right\vert \right)  ^{2}$ reduces finally to $1+\left\vert
\ln h\right\vert $ due to the inequality $\left\vert a\left(  \ln h\right)
\right\vert \leq c\left(  1+\left\vert \ln h\right\vert \right)
^{-1}\mathcal{N}$ inherited from (\ref{4.17}), (\ref{4.13}) and (\ref{4.19}).

The third term $I_{3}$ is examined in the simplest way. Indeed, formulas
(\ref{3.53}) and (\ref{3.11}), (\ref{3.77}) assure that%
\begin{align}
&  h^{-1/2}D\left(  \nabla\right)  \mathbf{W}\left(  h,x\right) \label{4.32}\\
&  =h^{1/2}(\left(  D_{\zeta}\mathcal{X}\left(  \zeta\right)  +\mathcal{Y}%
\left(  \zeta\right)  \right)  \mathcal{D}\left(  \nabla_{y}\right)
\mathbf{w}\left(  h,x\right)  +hD_{y}W^{2}\left(  \zeta,\nabla_{y}\right)
\mathbf{w}\left(  h,x\right)  )\nonumber
\end{align}
and we readily obtain that
\[
\left\vert I_{3}\right\vert \leq ch^{1/2}h^{1/2}\left\vert \ln h\right\vert
\mathcal{N}\mathbf{n}.
\]

Estimating the fourth term $I_{4}$ in (\ref{4.27}) follows the standard scheme
in the theory of thin plates with a modification caused by the cut-off
function $X_{h}^{\theta}$ present in (\ref{4.21}). The function $\widetilde
{w}^{\theta}$ satisfies the system%
\[
\mathcal{L}\left(  \nabla_{y}\right)  \widetilde{w}^{\theta}=\widetilde
{g}^{\theta}:=X_{h}^{\theta}g-[\mathcal{L}(\nabla_{y}),X_{h}^{\theta}]
\]
where $[\mathcal{L},X_{h}^{\theta}]$ is the commutator of the differential
operator with the cut-off function $X_{h}^{\theta}.$ By (\ref{4.14}) and
(\ref{4.25}), we have%
\begin{align}
&  \int_{\omega}\left(  S_{1,h}^{2}|(1-X_{h}^{\theta})g^{\prime}|^{2}%
+h^{-2}S_{2,h}^{2}|(1-X_{h}^{\theta})g_{3}|^{2}\right)  dy\leq ch^{2}\left(
1+\left\vert \ln h\right\vert \right)  ^{2}\mathcal{N},\label{4.35}\\
&  \int_{\omega}S_{1,h}^{2}|[L^{\prime},X_{h}^{\theta}]\widetilde{w}^{\prime
}|^{2}dy\leq c\int_{\mathbb{B}_{2R_{\theta}h}}\left(  r+h\right)  ^{2}\left(
1+\left\vert \ln\left(  r^{2}+h^{2}\right)  \right\vert \right)
^{2}\big(h^{-4}|\widetilde{w}^{\prime}|^{2}+h^{-2}|\nabla_{y}\widetilde
{w}^{\prime}|^{2}\big)dy\nonumber\\
&  \leq ch^{2}\left\vert \ln h\right\vert ^{4}\mathcal{N},\nonumber\\
&  \int_{\omega}S_{2,h}^{2}|[L_{3},X_{h}^{\theta}]\widetilde{w}_{3}%
|^{2}dy\nonumber\\
&  \leq c\int_{\mathbb{B}_{2R_{\theta}h}}\left(  r+h\right)  ^{4}\left(
1+\left\vert \ln\left(  r^{2}+h^{2}\right)  \right\vert \right)
^{2}\big(h^{-8}|\widetilde{w}_{3}|^{2}+h^{-6}|\nabla_{y}\widetilde{w}_{3}%
|^{2}+h^{-4}|\nabla_{y}^{2}\widetilde{w}_{3}|^{2}+h^{-2}|\nabla_{y}%
^{3}\widetilde{w}_{3}|^{2}\big)dy\nonumber\\
&  \leq ch^{2}\left\vert \ln h\right\vert ^{4}\mathcal{N}.\nonumber
\end{align}
In these calculations we used that weights are small on the disk
$\mathbb{B}_{2R_{\theta}h}$\ and coefficients of order $k$ derivatives in the
first- and third-order operators $[L^{\prime},X_{h}^{\theta}]$ and
$[L_{3},X_{h}^{\theta}]$ are $O(h^{k-2})$ and $O(h^{k-4}),$ respectively.
Returning to the notation in Section \ref{sect3.2}, we set $\widetilde{U}%
^{p}\left(  y,\zeta\right)  =W^{p}\left(  \zeta,\nabla_{y}\right)
\widetilde{w}^{\theta}\left(  y\right)  ,$ $p=0,1,2,$ and, similarly to
(\ref{4.32}), write
\begin{align}
I_{4}  &  =h^{-5/2}\big(AD_{\zeta}\widetilde{U}^{0},D\left(  \nabla\right)
\mathbf{v}^{\omega}\big)_{\Omega_{h}}+h^{-3/2}\big(A\big(D_{\zeta}%
\widetilde{U}^{1}+D_{y}\widetilde{U}^{0}\big),D\left(  \nabla\right)
\mathbf{v}^{\theta}\big)_{\Omega_{h}}\label{4.36}\\
&  +h^{-1/2}\big(A\big(D_{\zeta}\widetilde{U}^{2}+D_{y}\widetilde{U}%
^{1}\big),D\left(  \nabla\right)  \mathbf{v}^{\theta}\big)_{\Omega_{h}%
}+h^{1/2}\big(AD_{y}\widetilde{U}^{2},D\left(  \nabla\right)  \mathbf{v}%
^{\theta}\big)_{\Omega_{h}}\nonumber\\
&  =h^{-3/2}\left(  A\left(  D_{\zeta}\mathcal{X}+\mathcal{Y}\right)
\mathcal{D}\left(  \nabla_{y}\right)  w^{\theta},D_{\zeta}\mathbf{v}^{\omega
}\right)  _{\Omega_{h}}+h^{-1/2}\big(A\big(D_{\zeta}\widetilde{U}^{2}%
+D_{y}\widetilde{U}^{1}\big),D_{y}\mathbf{v}^{\omega}\big)_{\Omega_{h}%
}\nonumber\\
&  +\big(AD_{y}\widetilde{U}^{2},D_{\zeta}\mathbf{v}^{\omega}\big)_{\Omega
_{h}}+h^{1/2}\big(AD_{y}\widetilde{U}^{2},D_{y}\mathbf{v}^{\omega
}\big)_{\Omega_{h}}\nonumber\\
&  =:h^{-3/2}I_{40}+h^{-1/2}I_{41}+h^{1/2}I_{42}.\nonumber
\end{align}
Here, we have taken into account that $D_{\zeta}\widetilde{U}^{0}=0,$
$D_{\zeta}\widetilde{U}^{1}+D_{y}\widetilde{U}^{0}=0$ due to (\ref{3.8}) and,
moreover, $I_{40}=0$ because the $3\times6$-matrix $\mathcal{X}$\ solves
problem (\ref{3.12}) and $w^{\theta}$ is independent of $\zeta.$ Recalling
problem (\ref{3.99}) with $p=3,$ we see that%
\begin{align}
I_{41}  &  =-\big(D_{y}^{\top}A\big(D_{\zeta}\widetilde{U}^{2}+D_{y}%
\widetilde{U}^{1}\big)+D_{\zeta}^{\top}AD_{y}\widetilde{U}^{1},\mathbf{v}%
^{\theta}\big)_{\Omega_{h}}+%
{\textstyle\sum\nolimits_{\pm}}
\pm\big(D_{3}^{\top}AD_{y}\widetilde{U}^{2},\mathbf{v}^{\omega}\big)_{\Sigma
_{h}^{\pm}}\label{4.37}\\
&  =\big(D_{\zeta}^{\top}AD_{\zeta}\widetilde{U}^{3}+f^{0},\mathbf{v}^{\omega
}\big)_{\Omega_{h}}-%
{\textstyle\sum\nolimits_{\pm}}
\pm\big(D_{3}^{\top}AD_{\zeta}\widetilde{U}^{3},\mathbf{v}^{\omega
}\big)_{\Sigma_{h}^{\pm}}\nonumber\\
&  =\big(f^{0},\mathbf{v}^{\omega}\big)_{\Omega_{h}}-\big(AD_{\zeta}%
\widetilde{U}^{3},D_{\zeta}\mathbf{v}^{\omega}\big)_{\Omega_{h}}\nonumber\\
&  =\big(f^{0},\mathbf{v}^{\omega}\big)_{\Omega_{h}}-h\big(AD_{\zeta
}\widetilde{U}^{3},D\left(  \nabla\right)  \mathbf{v}^{\omega}\big)_{\Omega
_{h}}+h\big(AD_{\zeta}\widetilde{U}^{3},D_{y}\mathbf{v}^{\omega}%
\big)_{\Omega_{h}}.\nonumber
\end{align}
We emphasize that the differential properties (\ref{4.000}) put all entries of
scalar products into the Lebesgue space $L^{2}\left(  \Omega_{h}\right)  .$
Furthermore, we have the estimates
\begin{align}
h^{-1/2}\big|\left(  f^{0},\mathbf{v}^{\omega}\right)  _{\Omega_{h}}-\left(
f^{0},\mathbf{v}\right)  _{\Omega_{h}}\big|  &  =h^{-1/2}\left\vert
\sum\limits_{i=1}^{2}\left(  f_{i}^{0},\left(  1-X_{h}^{\omega}\right)
\mathbf{v}_{i}\right)  _{\Omega_{h}}\right\vert \label{4.38}\\
&  \leq ch^{-1/2}h^{1/2}\left\Vert f^{0};L^{2}\left(  \Omega_{1}\right)
\right\Vert h~|||\mathbf{v};\omega\left(  h\right)  \times\left(
-h/2,h/2\right)  |||_{\bullet}\nonumber\\
&  \leq ch\mathcal{N}\mathbf{n},\nonumber\\
h^{-1/2}\big|h\big(AD_{\zeta}\widetilde{U}^{3},D\left(  \nabla\right)
\mathbf{v}^{\omega}\big)_{\Omega_{h}}\big|  &  \leq ch^{1/2}h^{1/2}%
\mathcal{N}\left\Vert D\left(  \nabla\right)  \mathbf{v}^{\omega};L^{2}\left(
\Omega_{h}\right)  \right\Vert \leq ch\mathcal{N}\mathbf{n}.\nonumber
\end{align}
Notice that the factors $h^{1/2}$ came into (\ref{4.38}) due to the change
$z\mapsto\zeta=h^{-1}z.$ It remains to consider the last terms in (\ref{4.36})
and (\ref{4.37}). We write%
\begin{gather*}
h^{1/2}\big(A\big(D_{y}\widetilde{U}^{2}+D_{\zeta}\widetilde{U}^{3}%
\big),D_{y}\mathbf{v}^{\omega}\big)_{\Omega_{h}}=h^{1/2}\big(A\big(D_{y}%
\widetilde{U}^{2}+D_{\zeta}\widetilde{U}^{3}\big),D_{y}e_{3}\overline
{\mathbf{v}}_{3}^{\omega}\big)_{\Omega_{h}}\\
+h^{1/2}\big(A\big(D_{y}\widetilde{U}^{2}+D_{\zeta}\widetilde{U}%
^{3}\big),D_{y}\left(  \mathbf{v}^{\omega}-e_{3}\overline{\mathbf{v}}%
_{3}^{\omega}\right)  \big)_{\Omega_{h}}=:h^{1/2}I_{42}^{1}+h^{1/2}I_{42}^{2},
\end{gather*}
where $\overline{\mathbf{v}}_{3}^{\omega}\left(  y\right)  =h^{-1}%
{\textstyle\int_{-h/2}^{h/2}}
\mathbf{v}_{3}^{\omega}\left(  y,z\right)  dz.$ The known inequality%
\[
\left\Vert D_{y}\left(  \mathbf{v}^{\omega}-e_{3}\overline{\mathbf{v}}%
_{3}^{\omega}\right)  ;L^{2}\left(  \Omega_{h}\right)  \right\Vert
\leq\left\Vert D\left(  \nabla\right)  \mathbf{v}^{\omega};L^{2}\left(
\Omega_{h}\right)  \right\Vert ,\ \ \ \ \forall\mathbf{v}^{\omega}\in
H_{0}^{1}\left(  \Omega_{h},\upsilon_{h}\right)  ^{3},
\]
see, e.g., \cite[Section 4.3]{na273} and \cite[Proposition 3.3.13]{Nabook},
yields:
\[
h^{1/2}\left\vert I_{42}^{2}\right\vert \leq ch^{1/2}h^{1/2}\mathcal{N}%
\mathbf{n}.
\]

Finally, calculation (\ref{3.21}) together with formulas (\ref{3.16}),
(\ref{3.25}) for the matrix $\mathcal{A}$ and the relation $\mathcal{D}%
_{3}\left(  \nabla_{y}\right)  =2^{-1/2}\nabla_{y}\mathcal{D}^{\prime}\left(
\nabla_{y}\right)  $ between blocks in the matrix $\mathcal{D}\left(
\nabla_{y}\right)  ,$ see (\ref{3.10}), demonstrate that%
\begin{align*}
h^{1/2}I_{42}^{1}-h^{1/2}\left(  f_{3}^{1},\overline{\mathbf{v}}_{3}^{\omega
}\right)  _{\Omega_{h}}  &  =h^{3/2}\Bigg(2^{-1/2}\left(  \mathcal{D}^{\prime
}\left(  -\nabla_{y}\right)  ^{\top}\mathcal{A}_{\left(  3\right)
}\mathcal{D}_{3}\left(  \nabla_{y}\right)  \widetilde{w}_{3}^{\theta}%
,\nabla_{y}\overline{\mathbf{v}}_{3}^{\omega}\right)  _{\omega}\\
&  -\left(  f_{30}^{1},\overline{\mathbf{v}}_{3}^{\omega}\right)  _{\omega
}+\sum\limits_{i=1}^{2}\left(  f_{3i}^{1}+X_{h}^{\theta}%
{\displaystyle\int_{-1/2}^{1/2}}
\zeta f_{i}^{0}\left(  \cdot,\zeta\right)  d\zeta,\frac{\partial
\overline{\mathbf{v}}_{3}^{\omega}}{\partial y_{i}}\right)  _{\omega}\Bigg).
\end{align*}
We here took into account that in the approximate solution (\ref{4.21}) the
component $\widetilde{w}$ of (\ref{4.12}) is multiplied with the cut-off
function $X_{h}^{\theta},$ see (\ref{4.20}). Hence, to use the integral
identity (\ref{3.33}) with particular test functions $v_{i}=0$ and
$v_{3}=X_{h}^{\theta}\overline{\mathbf{v}}_{3}^{\omega},$ which vanish in
neighborhood of $\partial\omega$ and $\mathcal{O}$, we need to process several
commutators with coefficient supports located in the small disk $\overline
{\mathbb{B}_{2R_{\theta}h}}$ in the same way as we had done in (\ref{4.35}).
All of them receive bounds of similar order in $h,$ always much better than in
(\ref{4.30}), and therefore we list only typical estimates%
\begin{align*}
&  h^{3/2}\left\vert \big(\big[\mathcal{D}^{\prime}\left(  \nabla_{y}\right)
^{\top}\mathcal{A}_{\left(  3\right)  }\mathcal{D}_{3}\left(  \nabla
_{y}\right)  ,X_{h}^{\theta}\big]\widetilde{w}_{3},\nabla_{y}\overline
{\mathbf{v}}_{3}^{\omega}\big)_{\omega}\right\vert \\
&  \leq ch^{3/2}\left(  \int_{\mathbb{B}_{2R_{\theta}h}}\big(h^{-6}\left\vert
\widetilde{w}_{3}\right\vert ^{2}+h^{-4}\left\vert \nabla_{y}\widetilde{w}%
_{3}\right\vert ^{2}+h^{-2}\left\vert \nabla_{y}^{2}\widetilde{w}%
_{3}\right\vert ^{2}\big)dy\right)  ^{1/2}\left\Vert \nabla_{y}\overline
{\mathbf{v}}_{3}^{\omega};L^{2}\left(  \mathbb{B}_{2R_{\theta}h}\right)
\right\Vert \\
&  \leq ch^{3/2}\left\vert \ln h\right\vert \mathcal{N}h^{1/2}\left\Vert
\nabla_{y}\mathbf{v}_{3}^{\omega};L^{2}\left(  \mathbb{B}_{2R_{\theta}h}%
\times\left(  -h/2,h/2\right)  \right)  \right\Vert \leq ch\left\vert \ln
h\right\vert ^{2}\mathcal{N}\mathbf{n},\\
&  h^{3/2}\left\vert \left(  f_{30}^{1},\overline{\mathbf{v}}_{3}^{\omega
}\right)  _{\omega}-\left(  f_{30}^{1},X_{h}^{\theta}\overline{\mathbf{v}}%
_{3}^{\omega}\right)  _{\omega}\right\vert \\
&  \leq ch^{3/2}\left\Vert f_{30}^{1};L^{2}\left(  \mathbb{\omega}\right)
\right\Vert h^{1/2}\left\Vert \mathbf{v}_{3}^{\omega};L^{2}\left(
\mathbb{B}_{2R_{\theta}h}\times\left(  -h/2,h/2\right)  \right)  \right\Vert
\leq ch^{2}\left\vert \ln h\right\vert ^{2}\mathcal{N}\mathbf{n},\\
&  h^{3/2}\left\vert \left(  f_{3i}^{1},\overline{\mathbf{v}}_{3}^{\omega
}\partial_{i}X_{h}^{\theta}\right)  _{\omega}\right\vert \\
&  \leq ch^{3/2}\left\Vert f_{3i}^{1};L^{2}\left(  \mathbb{\omega}\right)
\right\Vert h^{-3/2}\left\Vert \mathbf{v}_{3}^{\omega};L^{2}\left(
\mathbb{B}_{2R_{\theta}h}\times\left(  -h/2,h/2\right)  \right)  \right\Vert
\leq ch\left\vert \ln h\right\vert ^{2}\mathcal{N}\mathbf{n}.
\end{align*}
In this way we finally obtain:
\[
h^{1/2}\big|I_{41}^{1}-\left(  f_{3}^{1},\overline{\mathbf{v}}_{3}^{\omega
}\right)  _{\Omega_{h}}\big|\leq ch\left\vert \ln h\right\vert ^{2}%
\mathcal{N}\mathbf{n}.
\]

We are in position to reckon calculations performed. Comparing bounds in all
estimates for terms on the right of (\ref{4.27}) we detect the worst one in
(\ref{4.30}), namely $ch^{1/2}\left\vert \ln h\right\vert \mathcal{N}%
\mathbf{n}.$ Recalling also supposition (\ref{FF5}), we finally derive from
(\ref{4.24}) and (\ref{4.27}), (\ref{3.1}) the formula%
\[
\left(  AD\left(  \nabla\right)  \mathbf{v},D\left(  \nabla\right)
\mathbf{v}\right)  _{\Omega_{h}}\leq ch^{1/2}\left\vert \ln h\right\vert
\big(\mathcal{N+}\widetilde{\mathcal{N}}\big)\left\Vert D\left(
\nabla\right)  \mathbf{v};L^{2}\left(  \Omega_{h}\right)  \right\Vert
\]
which together with the weighted anisotropic Korn inequality of Theorem 1 in
\cite{ButCaNa1} lead to the following assertion.

\begin{proposition}
\label{prop5.33}Under assumptions (\ref{3.2}) and (\ref{FF1}), (\ref{FF2}) on
the right-hand side (\ref{3.1}) in problem (\ref{B1.11})-(\ref{B1.14}) (or
(\ref{B1.16}) in the variational form) the difference $\mathbf{v}$ of the true
$u$ and approximate $\mathbf{u}$ solutions meets the estimate%
\begin{equation}
|||u-\mathbf{u};\Omega_{h}|||_{\bullet}\leq c\left\Vert D\left(
\nabla\right)  \left(  u-\mathbf{u}\right)  ;L^{2}\left(  \Omega_{h}\right)
\right\Vert \leq Ch^{1/2}\left\vert \ln h\right\vert \big(\mathcal{N+}%
\widetilde{\mathcal{N}}\big) \label{4.40}%
\end{equation}
where $\mathbf{u}$\ is determined in (\ref{4.21}), $\mathcal{N}$ and
$\widetilde{\mathcal{N}}$ in (\ref{FF4}) and (\ref{FF5}), $|||\cdot
|||_{\bullet}$ stands for the norm (\ref{B2.7}) and $C$ is a constant
independent of both, the right-hand side $f$ and the small parameter
$h\in\left(  0,h_{0}\right]  ,$ $h_{0}>0.$
\end{proposition}

\begin{remark}
\label{remBoLa} As verified in \cite{na109}, the error $O(h^{1/2})$ of the
Kirchhoff plate model (that is without the small support (\ref{B1.13})) cannot
be improved because of the boundary layer phenomenon near the edge. We will
see that the boundary layer near the support $\theta_{h}$ brings a
perturbation of much bigger order $O(\left\vert \ln h\right\vert ^{-1}).$
\end{remark}

\subsection{Theorem on asymptotics\label{sect4.5}}

The complicated form (\ref{4.21}) with various cut-off functions was
technically needed to satisfy the stable boundary conditions (\ref{B1.13}) and
(\ref{B1.14}) in the variational formulation (\ref{B1.16}) of the problem but
after verifying estimate (\ref{4.40}) we may get rid of some waste items while
keeping the same order of proximity. We consider two versions of the
simplified asymptotic structures, namely%
\begin{equation}
\mathbf{u}^{\theta}\left(  h,x\right)  =h^{-1/2}\mathcal{P}\left(
h^{-1}x\right)  a\left(  \ln h\right)  +h^{-3/2}\sum\limits_{p=0}^{2}%
h^{p}W^{p}\left(  \zeta,\nabla_{y}\right)  \widetilde{w}\left(  y\right)
\label{4.41}%
\end{equation}
with the notation in (\ref{4.21}) and%
\begin{equation}
\mathbf{u}^{\omega}\left(  h,x\right)  =h^{-3/2}\sum\limits_{p=0}^{2}%
h^{p}W^{p}\left(  \zeta,\nabla_{y}\right)  \big(\widehat{w}\left(  y\right)
+G^{\sharp\theta}\left(  h,x\right)  a\left(  \ln h\right)  \big)+h^{-1/2}%
\widetilde{\mathcal{P}}\left(  h^{-1}x\right)  a\left(  \ln h\right)
\label{4.42}%
\end{equation}
where, comparing formula (\ref{4.8}) with the content of Section 3.5 in
\cite{ButCaNa1}, we set%
\begin{align}
G^{\sharp\theta}\left(  h,x\right)   &  =(1-\chi(r/2hR_{\theta}))\Phi^{\sharp
}\left(  y\right)  +\widehat{G}^{\sharp}\left(  y\right)  ,\label{4.43}\\
\widetilde{\mathcal{P}}\left(  \xi\right)   &  =\mathcal{P}\left(  \xi\right)
-\sum\limits_{p=1}^{2}W^{p}\left(  \zeta,\nabla_{\eta}\right)  \left(
1-\chi\left(  \rho/2R_{\theta}\right)  \right)  \left(  \Phi^{\sharp}\left(
\eta\right)  +d^{\sharp}\left(  \eta,0\right)  C^{\sharp}\right)  .\nonumber
\end{align}

\begin{theorem}
\label{th4.51}Under conditions of Proposition \ref{prop5.33}, the vector
functions (\ref{4.41}), (\ref{4.42}) meet estimate (\ref{4.40}).
\end{theorem}

\textbf{Proof.} According to (\ref{4.43}), the approximate asymptotic
solutions (\ref{4.41}) and (\ref{4.42}) differ from each other only for the
singular term, cf. (\ref{4.3}),%
\[
\left(  1-\chi\left(  2r/hR_{\theta}\right)  \right)  \Phi^{\sharp}\left(
y\right)  =(1-\chi\left(  2\rho/R_{\theta}\right)  )H(\Phi^{\sharp}\left(
\eta\right)  +d^{\sharp}\left(  \eta,0\right)  \Psi\ln h),
\]
and the concomitant rearranging of the rigid motion
\[
h^{-1/2}H^{-1}d^{\sharp}(y,0)(C^{\sharp}-\Psi\ln h)a\left(  \ln h\right)  ,
\]
taking into account the matching condition (\ref{4.16}).

Hence, it is sufficient to check up the estimate for $\mathbf{u}^{\theta}.$
First of all, we observe that the discrepancy of (\ref{4.41}) in the Dirichlet
conditions (\ref{B1.14}) on the lateral side $\upsilon_{h}$ of the plate is
equal to%
\begin{equation}
h^{-3/2}\sum\limits_{p=0}^{2}h^{p}W^{p}\left(  \zeta,\nabla_{y}\right)
H\Upsilon^{\sharp}\left(  h^{-1}y\right)  +h^{-1/2}\widetilde{\widetilde
{\mathcal{P}}}\left(  \varepsilon^{-1}x\right)  +h^{1/2}W^{2}\left(
\zeta,\nabla_{y}\right)  \widetilde{w}\left(  y\right)  . \label{4.44}%
\end{equation}
Fast decay of the first two terms in (\ref{4.44}) and the small coefficient
$h^{1/2}$ in the third term allow us to withdraw from (\ref{4.21}) the cut-off
function $X_{h}^{\omega}$ as well as term (\ref{4.22}) which gets the small
factor $h$ in (\ref{4.21}) and vanishes in the disk $\mathbb{B}_{R_{\omega}%
/2}\ni y,$ that is at a distance from $\theta_{h}.$ This requires a standard
and simple calculation based on weights in norm (\ref{B2.7}) and has been
outlined in many publications (see, e.g., \cite{na273}, \cite[\S 4.3]{Nabook}
etc.). To take off from (\ref{4.21}) the other cut-off function $X_{h}%
^{\theta}$ which is null in the vicinity of $\theta_{h}$, becomes much more
delicate issue since the support of $\left\vert \nabla_{y}X_{h}^{\theta
}\right\vert $ is located very near singularities of the asymptotic terms. We
make use of the weighted estimate (\ref{4.14}) for $\widetilde{w}$ and observe
that the expressions $h^{p-3/2}W^{p}\left(  \zeta,\nabla_{y}\right)  \left(
1-X_{h}^{\theta}\left(  y\right)  \right)  \widetilde{w}\left(  y\right)  $
vanish outside the cylinder $\mathbb{B}_{2R_{\theta}h}\times\left(
-h/2,h/2\right)  $ where all weights in (\ref{4.14}) and (\ref{B2.7}) are big.
In this way we write%
\begin{align*}
&  h^{p-3/2}|||W^{p}(1-X_{h}^{\theta})\widetilde{w};\Omega_{h}|||_{\bullet}\\
&  \leq ch^{p-3/2}h^{1/2}\Bigg(\sum\limits_{q=0}^{p-1}h^{-\left(
p-1-q\right)  }h^{-1}\left\vert \ln h\right\vert ^{-1}||\nabla_{y}%
^{q}\widetilde{w}^{\prime};L^{2}(\mathbb{B}_{2R_{\theta}h})||\\
&  +\sum\limits_{q=0}^{p}h^{-\left(  p-q\right)  }h^{-1}|\ln h|^{-1}%
||\nabla_{y}^{q}\widetilde{w}_{3};L^{2}(\mathbb{B}_{2R_{\theta}h})||\Bigg)\\
&  \leq ch^{p-1}\Bigg(\sum\limits_{q=0}^{p-1}h^{-\left(  p-1-q\right)
}h^{1-q}||r^{-2+q}(1+\left\vert \ln r\right\vert )^{-1}\nabla_{y}%
^{q}\widetilde{w}^{\prime};L^{2}(\mathbb{B}_{2R_{\theta}h})||\\
&  +\sum\limits_{q=0}^{p}h^{-\left(  p-q\right)  }h^{2-q}||r^{-3+q}%
(1+\left\vert \ln r\right\vert )^{-1}\nabla_{y}^{q}\widetilde{w}_{3}%
;L^{2}(\mathbb{B}_{2R_{\theta}h})||\Bigg).
\end{align*}
Note that the factors $h^{1/2}$ and $h^{-1}|\ln h|^{-1}$ arrived due to
integration in $z\in\left(  -h/2,h/2\right)  $ and taking the weights $S_{h1}$
and $hS_{h2}$ into account. Plugging into norms the weights $r^{-2-q}\left(
1+\left\vert \ln r\right\vert \right)  ^{-1}$ and $r^{-3-q}\left(
1+\left\vert \ln r\right\vert \right)  ^{-1}$ according to (\ref{4.14}), we
have gained back $h^{2-q}\left\vert \ln h\right\vert $ and $h^{3-q}\left\vert
\ln h\right\vert ,$ respectively, that had led to the bound $chN$ which is
even better than in (\ref{4.40}). $\blacksquare$

\subsection{Further simplifications of asymptotic forms\label{sect4.6}}

By virtue of (\ref{4.14}), (\ref{4.17}) and the definition of the elastic
capacity potential in \cite{ButCaNa1}, the elastic energy of the boundary
layer term $h^{-1/2}\widetilde{\mathcal{P}}\left(  h^{-1}x\right)  a\left(
\ln h\right)  $ in (\ref{4.42}) gets order%
\begin{equation}
\big\Vert D\left(  \nabla\right)  h^{-1/2}\widetilde{\mathcal{P}}%
a;L^{2}(\Omega_{h})\big\Vert^{2}=O((h^{-1}h^{3}h^{-2}|a\left(  \ln h\right)
|^{2})=O(|\ln h|^{-2}) \label{4.45}%
\end{equation}
and therefore it cannot be omitted in the asymptotic expansion of the elastic
fields in the plate $\Omega_{h}$ when one attends to achieve an error estimate
with the power-law bound $c_{\delta}h^{\delta},$ $\delta>0.$ However,
contenting ourselves with the logarithmic precision $O(\left\vert \ln
h\right\vert ^{-1}),$ we may write much simpler asymptotics.

\begin{theorem}
\label{th4.52}Under conditions of Proposition \ref{prop5.33}, the regular
solution $\widehat{w}\in H^{1}\left(  \omega\right)  ^{2}\times H^{2}\left(
\omega\right)  $ of the two-dimensional Sobolev-Dirichlet problem
(\ref{3.15}),(\ref{3.23}), (\ref{3.24}) is related to the solution $u\in
H_{0}^{1}\left(  \Omega_{h};\Gamma_{h}\right)  ^{3}$ of the three-dimensional
problem (\ref{B1.11})-(\ref{B1.14}) by the inequality%
\[
|||u-h^{-3/2}\sum\limits_{p=0}^{2}h^{p}W^{p}\widehat{w};\Omega_{h}%
|||_{\bullet}\leq c\left\vert \ln h\right\vert ^{-1/2}(\mathcal{N+}%
\widetilde{\mathcal{N}}),
\]
where $W^{p}\left(  \zeta,\nabla_{y}\right)  $ are differential operators in
(\ref{3.53}) and $c$ is a constant independent of $h\in\left(  0,h_{0}\right]
$ and quantities $\mathcal{N},$ $\widetilde{\mathcal{N}}$ given in
(\ref{FF4}), (\ref{FF5}).
\end{theorem}

\textbf{Proof.} We still have to estimate the ingredients $h^{p-3/2}%
W^{p}\left(  \zeta,\nabla_{y}\right)  G^{\sharp\theta}\left(  h,x\right)
a\left(  \ln h\right)  ,\ p=0,1,2,$ of the asymptotic solution (\ref{4.22}).
We take into account representation (\ref{4.8}) of the Green matrix
$G^{\sharp}$ which is smooth in the punctured domain $\omega\setminus
\mathcal{O}$ and recall Remark \ref{remUW} about $W^{p}\left(  \zeta
,\nabla_{y}\right)  $ together with formulas (\ref{B346345}), (\ref{B3.49})
for the fundamental matrix. Then we write%
\begin{align}
h^{-3/2}|||G_{3}^{\sharp\theta}a;\Omega_{h}|||_{\bullet}  &  \leq
ch^{-3/2}h^{1/2}h\left(  1+\int_{\omega\setminus\mathbb{B}_{R_{\theta}h}%
}\left(  \frac{r^{2}\left(  1+\left\vert \ln r\right\vert \right)  ^{2}%
}{\left(  r^{2}+h^{2}\right)  ^{2}\left(  1+\left\vert \ln\left(  r+h\right)
\right\vert \right)  ^{2}}\right.  \right. \label{4.47}\\
&  \left.  \left.  +\frac{\left(  1+\left\vert \ln r\right\vert \right)  ^{2}%
}{\left(  r^{2}+h^{2}\right)  \left(  1+\left\vert \ln\left(  r+h\right)
\right\vert \right)  ^{2}}\right)  dy\right)  ^{1/2}\left\vert a\left(  \ln
h\right)  \right\vert \nonumber\\
&  \leq c\left\vert \ln h\right\vert ^{1/2}\left\vert \ln h\right\vert
^{-1}=c\left\vert \ln h\right\vert ^{-1/2},\nonumber\\
h^{-1/2}|||W^{1}G^{\sharp\theta}a;\Omega_{h}|||_{\bullet}  & \nonumber\\
\leq ch^{-1/2}h^{1/2}  &  \left(  1+\int_{\omega\setminus\mathbb{B}%
_{R_{\theta}h}}\frac{\left(  1+\left\vert \ln r\right\vert \right)  ^{2}%
}{\left(  r^{2}+h^{2}\right)  \left(  1+\left\vert \ln\left(  r+h\right)
\right\vert \right)  ^{2}}dy\right)  ^{1/2}\left\vert a\left(  \ln h\right)
\right\vert \leq c\left\vert \ln h\right\vert ^{-1/2},\nonumber\\
h^{1/2}|||W^{2}G^{\sharp\theta}a;\Omega_{h}|||_{\bullet}  & \nonumber\\
\leq ch^{1/2}h^{1/2}  &  \left(  1+\int_{\omega\setminus\mathbb{B}_{R_{\theta
}h}}\frac{\left(  1+\left\vert \ln r\right\vert \right)  ^{2}}{r^{2}\left(
r^{2}+h^{2}\right)  \left(  1+\left\vert \ln\left(  r+h\right)  \right\vert
\right)  ^{2}}dy\right)  ^{1/2}\left\vert a\left(  \ln h\right)  \right\vert
\leq c\left\vert \ln h\right\vert ^{-1/2}.\nonumber
\end{align}
For $p=0,$ we considered a specific input of the third component $u_{3}$ into
the weighted norm (\ref{B2.7}) while in the other cases $p=1,2$ the estimates
were derived in the standard way: $h^{1/2}$ comes from integration in $z$ and
$1$ is inserted to bound the regular part of $G^{\sharp}.$ Note that the first
two integrals in (\ref{4.47}) are $O\left(  \left\vert \ln h\right\vert
\right)  $ while the third one is $O\left(  h^{-1}\left\vert \ln h\right\vert
\right)  $ so that the relation $\left\vert a\left(  \ln h\right)  \right\vert
=O(\left\vert \ln h\right\vert ^{-1})$ is very important. $\blacksquare$

\bigskip

A similar effect of the catastrophic precision drop due to the Dirichlet
condition at a small region was underlined in \cite{na472} for an
homogenization problem in a perforated domain, too.

Our last simplification of asymptotic formulas is based on the following two
observations. First,%
\[
||h^{-1/2}\widetilde{\mathcal{P}}a;L^{2}\left(  \Omega_{h}\right)
||^{2}=O((h^{-1}h^{3}|a\left(  \ln h\right)  |^{2})=O(h^{2}|\ln h|^{-2})
\]
because the factors $h^{-2}$ and $h^{3}$ were introduced into (\ref{4.45}) due
to differentiation and integration in the fast variables $\xi=h^{-1}x.$
Second, logarithmic singularities do not lead the Green matrix $G^{\sharp}$
out from the spaces $L^{2}\left(  \omega\right)  $ and $L^{2}\left(
\Omega_{h}\right)  .$ In other words, the boundary layer may be excluded from
the asymptotic form but $G^{\sharp}$ may be included without any cut-off
function when dealing with the Lebesgue norms of displacements (for stresses
and strains this does not hold true, of course).

Let us formulate a simple but substantial consequence of our previous results.

\begin{theorem}
\label{th4.53}Under conditions of Proposition \ref{prop5.33}, the inequality%
\[
\sum\limits_{i=0}^{2}\left\Vert u_{i}-h^{-1/2}\left(  w_{i}-h^{-1}%
z\frac{\partial w_{3}}{\partial y_{i}}\right)  ;L^{2}\left(  \Omega
_{h}\right)  \right\Vert +h\left\Vert u_{3}-h^{-3/2}w_{3};L^{2}\left(
\Omega_{h}\right)  \right\Vert \leq ch^{1/2}\left\vert \ln h\right\vert
(\mathcal{N+}\widetilde{\mathcal{N}})
\]
is valid where $w\left(  \ln h;x\right)  $ is the singular solution
(\ref{4.11}) of the two-dimensional problem (\ref{3.15}), (\ref{3.23}),
(\ref{3.24}) with the column $a\left(  \ln h\right)  \in\mathbb{R}^{4}$
calculated in (\ref{4.17}) and $c$ is a constant independent of $h\in\left(
0,h_{0}\right]  $ and quantities $\mathcal{N},$ $\widetilde{\mathcal{N}}$
given in (\ref{FF4}), (\ref{FF5}).
\end{theorem}

\section{A variational asymptotic model of a plate with small
supports\label{sect5}}

\subsection{A symmetric unbounded operator and its adjoint\label{sect5.1}}

We intend to apply a traditional scheme to model media with small defect in
the framework of self-adjoint extension of differential operators, see the
pioneering paper \cite{BeFa} together with the review paper \cite{Pav} where a
physical terminology refers to "potentials of zero radii". Thanks to the
block-diagonal structure of the matrix $\mathcal{L}\left(  \nabla_{y}\right)
,$ cf. Lemma \ref{lem3.1}, all preparatory results below are substantiated by
the papers \cite{KaPa} and \cite{na239} where a scalar fourth-order
differential operator and a second-order system are considered. At the same
time, a straight way verification on them is an accessible task, too.

Let $\mathfrak{A}^{0}$ be an unbounded operator in the Hilbert space
$L^{2}\left(  \omega\right)  ^{3}$ with the differential expression
$\mathcal{L}\left(  \nabla_{y}\right)  $ and the domain is%
\begin{equation}
\mathfrak{D}\left(  \mathfrak{A}^{0}\right)  =\big\{w\in C_{c}^{\infty}\left(
\overline{\omega}\setminus\mathcal{O}\right)  ^{3}:\text{(\ref{3.23}) is
fulfilled}\big\}. \label{5.1}%
\end{equation}
Since the matrix $\mathcal{L}\left(  \nabla_{y}\right)  $ is elliptic and
block-diagonal, the closure $\mathfrak{A}=\overline{\mathfrak{A}^{0}}$ gets
the same differential expression but the domain%
\begin{align}
\mathfrak{D}\left(  \mathfrak{A}\right)   &  =\big\{w\in H^{2}\left(
\omega\right)  ^{2}\times H^{4}\left(  \omega\right)  :w\left(  y\right)
=0,\ \partial_{n}w_{3}\left(  y\right)  =0,\ y\in\partial\omega,\label{5.2}\\
&  \qquad\qquad w\left(  \mathcal{O}\right)  =0\in\mathbb{R}^{3},\ \nabla
_{y}w_{3}\left(  \mathcal{O}\right)  =0\in\mathbb{R}^{2}\big\}.\nonumber
\end{align}
Here, the Sobolev theorem on embedding $H^{2}\subset C$ is taken into account
as well as the fact that $H^{1}\left(  \omega\right)  \nsubseteq C\left(
\omega\right)  .$ In this way the property of $w$ in (\ref{5.1}) to vanish in
a neighborhood of the coordinate origin $\mathcal{O}$\ were converted into
five point conditions in (\ref{5.2}) including the Sobolev one (\ref{3.24}).

Both $\mathfrak{A}^{0}$\ and $\mathfrak{A}$\ are symmetric but not
self-adjoint. Indeed, their adjoint operator $\mathfrak{A}^{\ast}$ keeps, as
it is shown in \cite{KaPa, na239} and can be verified directly, the
differential expression $\mathcal{L}\left(  \nabla_{y}\right)  $ but its
domain becomes much bigger, namely%
\begin{align}
\mathfrak{D}\left(  \mathfrak{A}^{\ast}\right)   &  =\left\{
w=w_{\text{(sing)}}+w_{\text{(reg)}}:w_{\text{(reg)}}\in H^{2}\left(
\omega\right)  ^{2}\times H^{4}\left(  \omega\right)  ,\ w_{\text{(reg)}%
}\left(  y\right)  =0,\ \right. \label{5.3}\\
&  \qquad\qquad\partial_{n}w_{\text{(reg)}3}\left(  y\right)  =0,\ y\in
\partial\omega,\ w_{\text{(sing)}}^{\prime}\left(  y\right)  =G^{\prime
}\left(  y,\mathcal{O}\right)  a^{\prime},\ a^{\prime}=\left(  a_{1}%
,a_{2}\right)  ^{\top},\nonumber\\
&  \qquad\qquad w_{\text{(sing)}3}^{\prime}\left(  y\right)  =G_{3}\left(
y,\mathcal{O}\right)  a_{0}+G_{3}^{1}\left(  y,\mathcal{O}\right)  a_{3}%
+G_{3}^{2}\left(  y,\mathcal{O}\right)  a_{4},\nonumber\\
&  \qquad\qquad\left.  a_{0}\in\mathbb{R},\ a=\left(  a_{1},a_{2},a_{3}%
,a_{4}\right)  ^{\top}\in\mathbb{R}^{4}\right\}  .\nonumber
\end{align}
Recall that the Green matrix $G^{\prime}$ of problem (\ref{4.5}) and the Green
function $G_{3}$ of the Dirichlet problem for the operator $\mathcal{L}%
_{3}\left(  \nabla_{y}\right)  $ together with its derivatives (\ref{4.00})
live in the Lebesgue space but higher-order derivatives do not. In particular,
$\mathfrak{D}\left(  \mathfrak{A}^{\ast}\right)  \subset L^{2}\left(
\omega\right)  ^{3}.$

\subsection{Self-adjoint extensions\label{sect5.2}}

One readily observes that codimension of the subspace $\mathfrak{D}\left(
\mathfrak{A}\right)  \subset\mathfrak{D}\left(  \mathfrak{A}^{\ast}\right)  $
is $10=5+5,$ i.e. five free constants $a_{0},...,a_{4}$ in (\ref{5.3}) and
$b_{0}=w_{\text{(reg)}3}\left(  \mathcal{O}\right)  ,\ b=d^{\sharp}\left(
\nabla_{y},0\right)  w_{\text{(reg)}}\left(  \mathcal{O}\right)  \in
\mathbb{R}^{4}.$ Moreover, the defect index of $\mathfrak{A}^{\ast}$\ equals
$(5:5)$, cf. \cite{KaPa, na239}. Hence, according to the von Neumann theorem,
see \cite[\S 4.4]{BiSo}, the operator $\mathfrak{A}$ admits self-adjoint extensions.

The Friedrichs extension (or the hard one, \cite[\S 10.3]{BiSo}) is obtained,
for example, as a restriction of $\mathfrak{A}^{\ast}$ onto the subspace
$\left\{  w\in\mathfrak{D}\left(  \mathfrak{A}^{\ast}\right)  :a_{0}%
=0,\ a=0\in\mathbb{R}^{4}\right\}  $ of codimension $5$ and corresponds to the
classical formulation in $H^{2}\left(  \omega\right)  ^{2}\times H^{4}\left(
\omega\right)  $ of problem (\ref{3.15}),(\ref{3.23}) with the right-hand side
$g\in L^{2}\left(  \omega\right)  ^{3}.$ The Dirichlet-Sobolev problem
(\ref{3.15}),(\ref{3.23}), (\ref{3.24}), again with $g\in L^{2}\left(
\omega\right)  ^{3},$ is associated with a self-adjoint extension possessing
the domain%
\[
\big\{w\in\mathfrak{D}\left(  \mathfrak{A}^{\ast}\right)  :a=0\in
\mathbb{R}^{4},\ b_{0}=w_{3}\left(  \mathcal{O}\right)  \text{ but }%
a_{0}\text{ is arbitrary}\big\}.
\]
However, to realize modeling of the three-dimensional elasticity problem
(\ref{B1.11})-(\ref{B1.14}) in $\Omega_{h}$ we ought to deal with a different
self-adjoint extension.

\begin{theorem}
\label{th5.1}Let $M$ be a symmetric (real) $4\times4$-matrix. The restriction
$\mathfrak{A}_{M}$ of $\mathfrak{A}^{\ast}$ onto the linear space%
\begin{align}
\mathfrak{D}\left(  \mathfrak{A}_{M}\right)  =\big\{w\in\mathfrak{D}\left(
\mathfrak{A}^{\ast}\right)  : d^{\sharp}\left(  \nabla_{y},0\right)  ^{\top
}\widehat{w}\left(  \mathcal{O}\right)  =  &  Ma, \ \widehat{w}_{3}\left(
\mathcal{O}\right)  =0,\ \label{5.5}\\
&  \widehat{w}\left(  y\right)  =w_{\text{(reg)}}\left(  y\right)
+a_{0}G\left(  y,\mathcal{O}\right)  \big\}\nonumber
\end{align}
is a self-adjoint extension of the operator $\mathfrak{A}$ with domain
(\ref{5.2}).
\end{theorem}

\textbf{Proof.} To verify that (\ref{5.5}) constitutes the domain of a
self-adjoint operator, we first of all mention that $\mathfrak{D}\left(
\mathfrak{A}\right)  \subset\mathfrak{D}\left(  \mathfrak{A}_{M}\right)  $ and
$\dim\left(  \mathfrak{D}\left(  \mathfrak{A}^{\ast}\right)  /\mathfrak{D}%
\left(  \mathfrak{A}_{M}\right)  \right)  =5$ because five linear restrictions
are imposed on the free coefficients in (\ref{5.3}), namely on $a,$
$w_{\text{(reg)}3}\left(  \mathcal{O}\right)  $ and $d^{\sharp}\left(
\nabla_{y},0\right)  \widehat{w}\left(  \mathcal{O}\right)  =(\widehat{w}%
_{1}\left(  \mathcal{O}\right)  ,\widehat{w}_{2}\left(  \mathcal{O}\right)
,\partial_{1}\widehat{w}_{3}\left(  \mathcal{O}\right)  ,\partial_{2}%
\widehat{w}_{3}\left(  \mathcal{O}\right)  )^{\top},$ see (\ref{B3.49}). Then
the symplectic form%
\begin{equation}
\mathfrak{q}\left(  w,v\right)  =\left(  \mathfrak{A}^{\ast}w,v\right)
_{\omega}+\left(  w,\mathfrak{A}^{\ast}v\right)  _{\omega} \label{5.6}%
\end{equation}
is properly defined on the direct product. At the same time, with a clear
reason, form (\ref{5.6}) vanishes in the case both $w$ and $v$ belong to the
domain of a self-adjoint extension of $\mathfrak{A}.$ Thus, if we confirm that
$\mathfrak{q}\left(  w,v\right)  =0,\ \ \ \forall w,v\in\mathfrak{D}\left(
\mathfrak{A}_{M}\right)  $ the theorem is concluded.

We take functions from (\ref{5.3}) with the attributes $a_{\left(  w\right)
},a_{\left(  v\right)  }\in\mathbb{R}^{4}$ and write
\begin{align*}
\mathfrak{q}\left(  w,v\right)   &  =\lim_{t\rightarrow+0}\big(\left(
\mathcal{L}w,v\right)  _{\omega\setminus\mathbb{B}_{t}}-\left(  w,\mathcal{L}%
v\right)  _{\omega\setminus\mathbb{B}_{t}}\big)\\
&  =\lim_{t\rightarrow+0}\big(\left(  \mathcal{N}^{\prime}w^{\prime}%
,v^{\prime}\right)  _{\mathbb{S}_{t}}+\left(  \mathcal{N}_{3}w_{3},\left(
1,\nabla_{y}\right)  v_{3}\right)  _{\mathbb{S}_{t}}-\left(  w^{\prime
},\mathcal{N}^{\prime}v^{\prime}\right)  _{\mathbb{S}_{t}}-\left(  \left(
1,\nabla_{y}\right)  w_{3},\mathcal{N}_{3}v_{3}\right)  _{\mathbb{S}_{t}}\big)
\end{align*}
where $\mathbb{B}_{t}$ and $\mathbb{S}_{t}$ are the disk and circle with
radius $t$ and the center at $y=0$ while the Neumann operators $\mathcal{N}%
^{\prime}$ and $\mathcal{N}_{3}$ are taken from the Green formula for the
differential operator matrix $\mathcal{L}=\operatorname*{diag}\left\{
\mathcal{L}^{\prime},\mathcal{L}_{3}\right\}  $ in the domain $\omega
\setminus\mathbb{B}_{t}$ with a small hole. Natural integration properties of
the fundamental matrix (\ref{B346345}), given, e.g., by formulas (3.36) and
(3.42) in Section 3.4 in \cite{ButCaNa1}, provide a direct calculation of the
last limit according to the representation in (\ref{5.3}) and we finally
obtain%
\begin{align}
\mathfrak{q}\left(  w,v\right)   &  =\left(  \mathcal{L}w,v\right)  _{\omega
}-\left(  w,\mathcal{L}v\right)  _{\omega}\label{5.8}\\
&  =(d^{\sharp}\left(  \nabla_{y},0\right)  ^{\top}\widehat{v}\left(
\mathcal{O}\right)  +\mathcal{G}^{\sharp}a_{\left(  v\right)  })^{\top
}a_{\left(  w\right)  }-a_{\left(  w\right)  }^{\top}\big(d^{\sharp}\left(
\nabla_{y},0\right)  ^{\top}\widehat{w}\left(  \mathcal{O}\right)
+\mathcal{G}^{\sharp}a_{\left(  w\right)  }\big).\nonumber
\end{align}
Since the matrices $M$ and $\mathcal{G}^{\sharp}$\ are symmetric, inserting
relationship given in (\ref{5.5}) shows that the right-hand side of
(\ref{5.8}) vanishes. $\blacksquare$

\bigskip

We call (\ref{5.8}) the generalized Green formula, cf. \cite{na239, na161},
\cite[Ch.6]{NaPl}; it holds for any $w,v\in\mathfrak{D}\left(  \mathfrak{A}%
^{\ast}\right)  .$

\subsection{The operator model for the supported plate\label{sect5.3}}

We simplify representation (\ref{3.1}) of the right-hand side in the
three-dimensional problem (\ref{B1.11})-(\ref{B1.14}) and set%
\begin{equation}
f\left(  h,x\right)  =h^{-1/2}\left(  f_{1}^{0}\left(  y\right)  ,f_{2}%
^{0}\left(  y\right)  ,hf_{3}^{0}\left(  y\right)  \right)  ^{\top
},\ \ \ f_{j}^{0}\in L^{2}\left(  \omega\right)  . \label{5.9}%
\end{equation}
The introduced restrictions (\ref{3.2}) and (\ref{FF1}), (\ref{FF2}) are, of
course, satisfied and, according to (\ref{3.17}), (\ref{3.22}), (\ref{FF4})%
\begin{equation}
g\left(  y\right)  =f^{0}\left(  y\right)  =\left(  f_{1}^{0}\left(  y\right)
,f_{2}^{0}\left(  y\right)  ,f_{3}^{0}\left(  y\right)  \right)  ^{\top
},\ \ \ \mathcal{N}=\left\Vert f^{0};L^{2}\left(  \omega\right)  \right\Vert
,\ \ \ \widetilde{\mathcal{N}}=0. \label{5.10}%
\end{equation}

Let $\mathfrak{A}\left(  \ln h\right)  $ be the self-adjoint operator
$\mathfrak{A}_{M^{\sharp}\left(  \ln h\right)  }$ given by Theorem \ref{th5.1}
with the numeral matrix (\ref{4.18}). We consider the abstract equation%
\begin{equation}
\mathfrak{A}\left(  \ln h\right)  w=g\in L^{2}\left(  \omega\right)  ^{3}.
\label{5.11}%
\end{equation}
A solution of this equation takes form (\ref{4.11}) where $\widehat{w}\in
H^{1}\left(  \omega\right)  ^{2}\times H^{2}\left(  \omega\right)  $ is the
unique solution of the Sobolev-Dirichlet problem (\ref{3.15}), (\ref{3.23}),
(\ref{3.24}) and the column $a=a\left(  \ln h\right)  =M^{\sharp}\left(  \ln
h\right)  ^{-1}d^{\sharp}\left(  \nabla_{y},0\right)  \widehat{w}\left(
\mathcal{O}\right)  $ is perfectly defined through formulas (\ref{4.17}) and
(\ref{4.12}) for a small $h>0.$ In other words, the abstract equation
(\ref{5.11}) is uniquely solvable and, moreover, its solution $w\in
\mathfrak{D}\left(  \mathfrak{A}\left(  \ln h\right)  \right)  $ coincides
with generator (\ref{4.11}) of the outer asymptotic expansion (\ref{3.3})
constructed in Sections \ref{sect3.1} and \ref{sect4.1}.

Theorem \ref{th4.53} assures the following result.

\begin{corollary}
\label{cor5.99}A solution $w$ of equation (\ref{5.11}) with the chosen
self-adjoint extension $\mathfrak{A}\left(  \ln h\right)  $ and a right-hand
side as in (\ref{5.9}) and (\ref{5.10}) satisfies the estimate%
\[
\int_{\Omega_{h}}\sum\limits_{i=1}^{2}\big|u_{i}\left(  h,x\right)
-h^{-1/2}w_{i}\left(  \ln h,y\right)  \big|^{2}+h^{2}\big|u_{3}\left(
h,x\right)  -h^{-3/2}w_{3}\left(  \ln h,y\right)  \big|^{2}dx\leq
ch^{1/2}\left\vert \ln h\right\vert \mathcal{N}%
\]
where $\mathcal{N}$\ is given in (\ref{5.10}) and $c$ is a constant
independent of data (\ref{5.9}) and the small parameter $\varepsilon
\in(0,\varepsilon_{0}),$ $\varepsilon_{0}>0.$
\end{corollary}

\subsection{The variational model for the supported plate\label{sect5.4}}

The self-adjoint operator $\mathfrak{A}\left(  \ln h\right)  $ traditionally
gives rise to the energy functional%
\begin{equation}
\tfrac{1}{2}\left(  \mathfrak{A}\left(  \ln h\right)  w,w\right)  _{\omega
}-\left(  g,w\right)  _{\omega}. \label{5.13}%
\end{equation}
We also introduce the functional%
\begin{equation}
\mathfrak{E}\left(  w;g\right)  =\tfrac{1}{2}\left(  \mathfrak{A}^{\ast
}w,w\right)  _{\omega}-\left(  g,w\right)  _{\omega}+\tfrac{1}{2}a^{\top
}(M^{\sharp}\left(  \ln h\right)  a-d^{\sharp}\left(  \nabla_{y},0\right)
\widehat{w}\left(  \mathcal{O}\right)  ) \label{5.14}%
\end{equation}
which according to (\ref{5.5}) at $M=M^{\sharp}\left(  \ln h\right)  ,$
coincides with (\ref{5.13}) for $v\in\mathfrak{D}\left(  \mathfrak{A}\left(
\ln h\right)  \right)  $ but is defined on the whole space%
\begin{equation}
\mathfrak{H}=\big\{v\in\mathfrak{D}\left(  \mathfrak{A}^{\ast}\right)
:\widehat{w}_{3}\left(  \mathcal{O}\right)  =0\big\} \label{5.15}%
\end{equation}
which becomes Hilbert with the norm%
\[
\left\Vert w;\mathfrak{H}\right\Vert =\big(\big\Vert w_{\text{(reg)}}%
;H^{2}\left(  \omega\right)  ^{2}\times H^{4}\left(  \omega\right)
\big\Vert^{2}+\left\vert a_{0}\right\vert ^{2}+\left\vert a\right\vert
^{2}\big)^{1/2}.
\]

\begin{theorem}
\label{th5.2}A solution $w\in\mathfrak{D}\left(  \mathfrak{A}\left(  \ln
h\right)  \right)  $ of equation (\ref{5.11}) is the only stationary point of
functional (\ref{5.14}) in $\mathfrak{H}$.
\end{theorem}

\textbf{Proof.} Since, by definition, $\mathfrak{A}^{\ast}G^{\sharp}=0\in
L^{2}\left(  \omega\right)  ^{4},$ the generalized Green formula (\ref{5.8})
and the symmetry of the matrices $M^{\sharp}\left(  \ln h\right)  $ and
$\mathcal{G}^{\sharp}$ convert the variation of $\mathfrak{E}$ into%
\begin{align}
&  \frac{1}{2}\left(  \mathcal{L}w,v\right)  _{\omega}+\frac{1}{2}\left(
\mathcal{L}v,w\right)  _{\omega}-\left(  g,v\right)  _{\omega}+\frac{1}%
{2}a_{\left(  v\right)  }^{\top}\left(  M^{\sharp}\left(  \ln h\right)
a_{\left(  w\right)  }-d^{\sharp}\left(  \nabla_{y},0\right)  \widehat
{w}\left(  \mathcal{O}\right)  \right) \label{5.17}\\
&  +\frac{1}{2}a_{\left(  w\right)  }^{\top}\left(  M^{\sharp}\left(  \ln
h\right)  a_{\left(  v\right)  }-d^{\sharp}\left(  \nabla_{y},0\right)
\widehat{v}\left(  \mathcal{O}\right)  \right) \nonumber\\
&  =\frac{1}{2}\left(  \mathcal{L}w,v\right)  _{\omega}+\frac{1}{2}\left(
v,\mathcal{L}w\right)  _{\omega}-\left(  g,v\right)  _{\omega}+\frac{1}%
{2}a_{\left(  w\right)  }^{\top}\left(  d^{\sharp}\left(  \nabla_{y},0\right)
\widehat{v}\left(  \mathcal{O}\right)  +\mathcal{G}^{\sharp}a_{\left(
v\right)  }\right) \nonumber\\
&  -\frac{1}{2}a_{\left(  v\right)  }^{\top}\left(  d^{\sharp}\left(
\nabla_{y},0\right)  \widehat{w}\left(  \mathcal{O}\right)  +\mathcal{G}%
^{\sharp}a_{\left(  w\right)  }\right)  +\frac{1}{2}a_{\left(  v\right)
}^{\top}\left(  M^{\sharp}\left(  \ln h\right)  a_{\left(  w\right)
}-d^{\sharp}\left(  \nabla_{y},0\right)  \widehat{w}\left(  \mathcal{O}%
\right)  \right) \nonumber\\
&  +\frac{1}{2}a_{\left(  w\right)  }^{\top}\left(  M^{\sharp}\left(  \ln
h\right)  a_{\left(  v\right)  }-d^{\sharp}\left(  \nabla_{y},0\right)
\widehat{v}\left(  \mathcal{O}\right)  \right) \nonumber\\
&  =\left(  \mathcal{L}\widehat{w}-g,v\right)  _{\Omega}+a_{\left(  v\right)
}^{\top}\left(  M^{\sharp}\left(  \ln h\right)  a_{\left(  w\right)
}-d^{\sharp}\left(  \nabla_{y},0\right)  \widehat{w}\left(  \mathcal{O}%
\right)  \right)  .\nonumber
\end{align}
Equating (\ref{5.17}) to null and taking a test function $v\in C_{c}^{\infty
}\left(  \omega\setminus\mathcal{O}\right)  ^{3}$ yields%
\[
(\mathcal{L}\widehat{w}-g,v)_{\Omega}=0\Longrightarrow\mathcal{L}\left(
\nabla_{y}\right)  \widehat{w}\left(  y\right)  =g\left(  y\right)
,\ \ \ y\in\omega\setminus\mathcal{O},
\]
while conditions (\ref{3.23}) and (\ref{3.24}) are inherited from the
inclusion $w\in\mathfrak{H}$. Since $a_{\left(  v\right)  }\in\mathbb{R}^{4}$
is arbitrary, we conclude that expression (\ref{5.17}) vanishes provided
$M^{\sharp}\left(  \ln h\right)  a_{\left(  w\right)  }=d^{\sharp}\left(
\nabla_{y},0\right)  \widehat{w}\left(  \mathcal{O}\right)  .$ Thus, $w$ falls
into the domain $\mathfrak{D}\left(  \mathfrak{A}\left(  \ln h\right)
\right)  .$

In a similar way one verifies that a solution $w$ of (\ref{5.11}) annuls the
variation of functional (\ref{5.14}). $\blacksquare$

\bigskip

From this theorem it follows that the abstract equation (\ref{5.11}) is
equivalent to the variational problem for the quadratic functional
(\ref{5.14}), cf. \cite{MikhQF}. The latter is much more suitable for, e.g.,
numerical implementation because the space (\ref{5.15}) does not involve any
linear restriction on $w$ except for the Sobolev condition (\ref{3.24}).

The potential energy of the three-dimensional plate (\ref{B1.1}) under the
volume force (\ref{5.9}) is equal to%
\[
E_{h}\left(  u;f\right)  =\tfrac{1}{2}\left(  AD\left(  \nabla_{x}\right)
u,D\left(  \nabla_{x}\right)  u\right)  _{\Omega_{h}}-\left(  u,f\right)
_{\Omega_{h}}=-\tfrac{1}{2}\left(  u,f\right)  _{\Omega_{h}};
\]
here the Green formula for a solution $u$ of (\ref{B1.11})-(\ref{B1.14}) was
applied. The last expression demonstrates that, although the introduced model
cannot provide "good" approximation of the strain and stress fields near the
support area $\theta_{h},$ Corollary \ref{cor5.99} leads to a simple
asymptotic formula for the energy.

\begin{corollary}
\label{cor5.77}For the right-hand sides (\ref{5.9}) and (\ref{5.10}), the
solutions $u\left(  h,x\right)  $ of the three-dimensional problem
(\ref{B1.11})-(\ref{B1.14}) and the solution $w\left(  \ln h,y\right)  $ of
its two-dimensional model provide the following relationship between the
corresponding energy functionals%
\[
\left\vert E_{h}\left(  u;f\right)  -\mathfrak{E}\left(  w;g\right)
\right\vert \leq c\sqrt{\varepsilon}\mathcal{N}%
\]
where $\mathcal{N}$\ is given in (\ref{5.10}) and $c$ is a constant
independent of data (\ref{5.9}) and the small parameter $\varepsilon\in\left(
0,\varepsilon_{0}\right)  ,$ $\varepsilon_{0}>0.$
\end{corollary}

Modifying calculation (\ref{5.17}) and applying the standard Green formula in
$\omega,$ we derive that%
\begin{align}
\mathfrak{E}\left(  w;g\right)   &  =\tfrac{1}{2}(\mathcal{L}\widehat
{w},\widehat{w}+G^{\sharp}a)_{\omega}-(g,\widehat{w}+G^{\sharp}a)_{\omega
}+0\label{5.20}\\
&  =\tfrac{1}{2}(\mathcal{L}\widehat{w},\widehat{w})_{\omega}-(g,\widehat
{w})_{\omega}-\tfrac{1}{2}a^{\top}d\left(  \nabla_{y},0\right)  \widehat
{w}\left(  \mathcal{O}\right)  +a^{\top}d\left(  \nabla_{y},0\right)
\widehat{w}\left(  \mathcal{O}\right) \nonumber\\
&  =\mathcal{E}(\widehat{w};g)+\tfrac{1}{2}a^{\top}M^{\sharp}\left(  \ln
h\right)  a,\nonumber
\end{align}
where%
\begin{equation}
\mathcal{E}(\widehat{w};g)=\frac{1}{2}(\mathcal{AD}\left(  \nabla_{y}\right)
\widehat{w},\mathcal{D}\left(  \nabla_{y}\right)  \widehat{w})_{\omega
}-(g,\widehat{w})_{\omega}. \label{5.21}%
\end{equation}

In other words, the energy functional (\ref{5.14}) computed for the solution
$w\in\mathfrak{D}\left(  \mathfrak{A}\left(  \ln h\right)  \right)  $ of the
two-dimensional plate model (\ref{5.11}) is equal to the sum of the potential
energy (\ref{5.21}), stored by the Kirchhoff plate with the Sobolev-Dirichlet
conditions, and the correction term%
\begin{equation}
\frac{1}{2}a^{\top}M^{\sharp}\left(  \ln h\right)  a \label{5.22}%
\end{equation}
which describes the energy concentrated in the vicinity of the small clamped
zone $\theta_{h}.$ It should be mentioned that, in view of (\ref{4.19}) and
(\ref{4.18}), value (\ref{5.22}) gets order $O(\left\vert \ln h\right\vert
^{-1})$ and is positive. The latter complies extension (\ref{B1.3}) of the
Dirichlet area in the minimization problem%
\[
E_{h}\left(  u;f\right)  =\min_{v\in H_{0}^{1}\left(  \Omega_{h};\Gamma
_{h}\right)  ^{3}}E_{h}\left(  v;f\right)
\]
compared with the traditional problem on the plate (\ref{B1.1}) clamped along
the lateral side $\upsilon_{h}.$

Although the model energy functional (\ref{5.20}) gains the positive increment
(\ref{5.22}), the stationary point indicated in Theorem \ref{th5.2} does not
constitute its minimum because of the calculation%
\begin{align*}
\mathfrak{E}\left(  w+v;g\right)  -\mathfrak{E}\left(  w;g\right)   &
=\frac{1}{2}\left(  \mathcal{L}\widehat{v},\widehat{v}+G^{\sharp}a_{\left(
v\right)  }\right)  _{\omega}+\frac{1}{2}a_{\left(  v\right)  }^{\top
}M^{\sharp}\left(  \ln h\right)  a_{\left(  v\right)  }\\
&  =\frac{1}{2}\left(  \mathcal{AD}\left(  \nabla_{y}\right)  \widehat
{v},\mathcal{D}\left(  \nabla_{y}\right)  \widehat{v}\right)  _{\omega}%
+\frac{1}{2}a_{\left(  v\right)  }^{\top}\left(  M^{\sharp}\left(  \ln
h\right)  a_{\left(  v\right)  }-d^{\sharp}\left(  \nabla_{y},0\right)
\widehat{v}\left(  \mathcal{O}\right)  \right)
\end{align*}
where the last terms can get any sign in $\mathfrak{H}$. This terms vanishes
in $\mathfrak{D}\left(  \mathfrak{A}\left(  \ln h\right)  \right)  $ and,
quite expected, the energy functional (\ref{5.14}) admits the global minimum
over the intrinsic linear space (\ref{5.5}) at the solution $w$ of equation
(\ref{5.11}).

\subsection{Conclusive remarks\label{sect5.5}}

Asymptotics derived and justified in \S \ref{sect4} indicates a distinguishing
feature of the influence of the small support $\theta_{h}$ on the
strain-stress state of the plate $\Omega_{h},$ see (\ref{B1.1}) and
(\ref{B1.2}). Namely, the three-dimensional boundary layer phenomenon in the
vicinity of $\theta_{h}$ overrides the standard error estimate $O(h^{1/2})$
for the two-dimensional Kirchhoff model and the convergence rate $O(\left\vert
\ln h\right\vert ^{-1/2})$ in%
\begin{align}
h^{1/2}u_{i}\left(  h,y,h\zeta\right)   &  \rightarrow\widehat{w}_{i}\left(
y\right)  -\zeta\frac{\partial\widehat{w}_{3}}{\partial y_{i}}%
,\ i=1,2,\ \label{5.31}\\
h^{3/2}u_{3}\left(  h,y,h\zeta\right)   &  \rightarrow\widehat{w}_{3}\left(
y\right)  \text{\ in }H^{1}\left(  \omega\times\left(  -1/2,1/2\right)
\right)  ,\nonumber
\end{align}
see Theorem \ref{th4.52}, becomes too sluggish so that cannot maintain an
engineering application. An acceptable error estimate in Theorem \ref{th4.51},
comparable with the usual one in the Kirchhoff theory, requires to include
into the asymptotic form the three-dimensional elastic field, that is the
elastic logarithmical capacity potential $\mathcal{P}\left(  \varepsilon
^{-1}x\right)  $ in Section 3.5 of \cite{ButCaNa1} which is not described
explicitly yet even for the disk $\theta_{h}$ of radius $h$ and the isotropic
material. However, replacing in (\ref{5.31}) the regular solution $\widehat
{w}\in H^{1}\left(  \omega\right)  ^{2}\times H^{2}\left(  \omega\right)  $ of
the Dirichlet-Sobolev problem (\ref{3.15}), (\ref{3.23}), (\ref{3.24}) for the
singular solution (\ref{4.11}) procures the convergence rate $O\left(
h^{1/2}\left\vert \ln h\right\vert \right)  $ in$\ $%
\begin{align}
h^{1/2}u_{i}\left(  h,y,h\zeta\right)  -w_{i}\left(  \ln h,y\right)
+\zeta\tfrac{\partial\widehat{w}_{3}}{\partial y_{i}}\left(  \ln h,y\right)
&  \rightarrow0,\ \ \ \ i=1,2,\label{5.32}\\
h^{3/2}u_{3}\left(  h,y,h\zeta\right)  -w_{3}\left(  \ln h,y\right)   &
\rightarrow0\text{\ in }L^{2}\left(  \omega\times\left(  -1/2,1/2\right)
\right)  ,\nonumber
\end{align}
however in the Lebesgue, not Sobolev norm.

The problem to determine the necessary singular solution $w\in L^{2}\left(
\omega\right)  ^{2}\times H^{1}\left(  \omega\right)  $ is well-posed with
different formulations in Sections \ref{sect5.3} and \ref{sect5.4}. Its
structure (\ref{4.11}), (\ref{4.17}) of $w,$ refers to the Green matrix
(\ref{4.8}), (\ref{4.9}) and the integral characteristics $C^{\sharp}\left(
A,\theta\right)  $ in the representation of the logarithmic elastic potential
$\mathcal{P}$\ in the elastic layer, see Theorem 13 in \cite{ButCaNa1}. The
elastic logarithmic capacity matrix is a numeral $4\times4$-matrix and can be
computed numerically.

Convergence (\ref{5.32}) does not provide information about the stress and
strain fields in the plate, however, using a technique of local estimates in
\cite{na224, Nabook}, it is possible to verify that, under certain
restrictions on the volume forces (\ref{FF1}), this convergence becomes
pointwise outside a neighborhood of $\overline{\theta_{h}}.$ We mention that,
as discovered in \cite{na109}, see also \cite{na333}, pointwise convergence of
stresses and strains is always broken near the clamped edge $\upsilon_{h}$ of
the plate with or without small perturbation at $\theta_{h}.$

Both, the asymptotic procedure and the models, can be easily adapted to the
plate (\ref{B1.1}) with several small support areas $\theta_{h}^{1}%
,...,\theta_{h}^{J}$ sparsely distributed on one or two bases $\Sigma_{h}%
^{\pm}.$

In view of block-diagonal structure (\ref{3.25}) of the matrix $\mathcal{A}$
the Dirichlet-Sobolev problem (\ref{3.15}), (\ref{3.23}), (\ref{3.24})
naturally decouples into independent problems for the vector $w^{\prime
}=(w_{1},w_{2})$ of longitudinal displacements and the deflection $w_{3}.$
Since the elasticity problem for the infinite layer $\Lambda=\mathbb{R}%
^{2}\times(-1/2,1/2)$ clamped over the only area $\theta\times\{-1/2\}$ loses
the symmetry with respect to the plane $\left\{  \xi\in\mathbb{R}^{3}:\xi
_{3}=\zeta=0\right\}  ,$ the elastic capacity matrix $C^{\sharp}$ does not get
in general a block-diagonal structure and thus the relationship imposed in
(\ref{5.5}) on the regular and similar components of the solution
$w=(w^{\prime},w_{3})$ of the model equation (\ref{5.11}) couples the vector
$w^{\prime}$ and the scalar $w_{3}$ at the level $O(\left\vert \ln
h\right\vert ^{-1}).$

\end{document}